\RequirePackage{fix-cm}

\documentclass[smallextended]{svjour3} 
\smartqed
\usepackage{graphicx}
\usepackage[utf8]{inputenc}
\usepackage{hyperref}
\hypersetup{colorlinks=true}
\usepackage{amsmath}
\usepackage{amsfonts}
\usepackage{amssymb}
\usepackage{subcaption}
\usepackage{listings}
\usepackage{graphicx}
\usepackage{appendix}
\usepackage{enumerate}
\usepackage{grffile}
\usepackage{url}
\usepackage{stfloats}

\usepackage[linesnumbered,ruled,vlined]{algorithm2e}
\usepackage{multirow}
\usepackage{multicol}

\begin{document}

\title{Optimal age-specific vaccination control for COVID-19: an Irish case study}
\author{Eleni Zavrakli        \and
        Andrew Parnell          \and
        David Malone            \and 
        Ken Duffy               \and 
        Subhrakanti Dey         
}


\institute{Eleni Zavrakli \at
              Hamilton Institute, Maynooth University, Co. Kildare, Ireland \\
             I-Form Advanced Manufacturing Research Centre, Ireland \\ 
            Department of Mathematics and Statistics, Maynooth University, Co. Kildare, Ireland \\
              \email{Eleni.Zavrakli@mu.ie}
           \and
           Andrew Parnell \at
              Hamilton Institute, Maynooth University, Co. Kildare, Ireland \\
             I-Form Advanced Manufacturing Research Centre, Ireland \\ 
            Department of Mathematics and Statistics, Maynooth University, Co. Kildare, Ireland 
            \and 
            David Malone \at 
            Hamilton Institute, Maynooth University, Co. Kildare, Ireland \\
            Department of Mathematics and Statistics, Maynooth University, Co. Kildare, Ireland 
            \and 
            Ken Duffy \at 
            Hamilton Institute, Maynooth University, Co. Kildare, Ireland 
            \and 
            Subhrakanti Dey \at 
            Department of Electrical Engineering, Uppsala University, Sweden\\
}

\date{}

\maketitle

\begin{abstract}
The outbreak of a novel coronavirus causing severe acute respiratory syndrome in December 2019 has escalated into a worldwide pandemic. In this work, we propose a compartmental model to describe the dynamics of transmission of infection and use it to obtain the optimal vaccination control. The model accounts for the various stages of the vaccination and the optimisation is focused on minimising the infections to protect the population and relieve the healthcare system. As a case study we selected the Republic of Ireland. We use data provided by Ireland's COVID-19 Data-Hub and simulate the evolution of the pandemic with and without the vaccination in place for two different scenarios, one representative of a national lockdown situation and the other indicating looser restrictions in place. One of the main findings of our work is that the optimal approach would involve a vaccination programme where the older population is vaccinated in larger numbers earlier while simultaneously part of the younger population also gets vaccinated to lower the risk of transmission between groups. We compare our simulated results with that of the vaccination policy taken by the Irish government to explore the advantages of our optimisation method. Our comparison suggests that a similar reduction in cases may have been possible even with a reduced set of vaccinations being available for use.
\end{abstract}


\maketitle

\section*{Declarations}
\subsection*{Funding}
This publication has emanated from research supported in part by a research grant from Science Foundation Ireland (SFI) under Grant Number 16/RC/3872 and is co-funded under the European Regional Development Fund. This work was also supported by a Science Foundation Ireland COVID-19 Rapid Response Grant number 20/COV/0081. Andrew Parnell’s work was additionally supported by: a Science Foundation Ireland Career Development Award (17/CDA/4695); an investigator award (16/IA/4520); a Marine Research Programme funded by the Irish Government, co-financed by the European Regional Development Fund (Grant-Aid Agreement No. PBA/CC/18/01); European Union’s Horizon 2020 research and innovation programme InnoVar under grant agreement No 818144; SFI Centre for Research Training in Foundations of Data Science 18CRT/6049, and SFI Research Centre award 12/RC/2289\_P2.

\subsection*{Conflicts of interest}
The authors have no conflicts of interest to declare that are relevant to the content of this article.

\subsection*{Code availability}
The code used for our simulations was written in R and can be found on:
\url{https://github.com/elenizavrakli/Optimal-Age-Specific-Vaccination-Control}

\vfill


\section{Introduction}
\label{intro}
In late 2019, an outbreak of pneumonia of unknown cause was reported in the city of Wuhan in the Hubei province of China \cite{cite-key,Who1,doi:10.1056/NEJMoa2001017,https://doi.org/10.1002/jmv.25678,cluster,newcovid,doi:10.1056/NEJMoa2001316}. The virus was named Severe Acute Respiratory Syndrome Coronavirus 2 (SARS-CoV-2) by the World Health Organisation \cite{Who2} and the disease that it causes is referred to as COVID-19 \cite{Who3}. The disease quickly became a source of international worry as it spread around China with most countries around the world following \cite{geogr}. By March 2020 most countries in the world had confirmed cases of COVID-19, including the Republic of Ireland, where the first confirmed case was on February 29th 2020 \cite{times}, the same day that WHO raised the risk warning for the virus to "very high" \cite{times2}. In the absence of an effective treatment or a vaccine, governments worldwide started implementing protective measures with most of them announcing a national lockdown to try and control the spread of the virus and reduce the strain on their healthcare systems. In Ireland, the first restrictive measures and social distancing guidelines were first announced on March 12th 2020 and were initially intended to last until March 29th. While countries tried to control the virus and protect their citizens, the scientific community committed to coming up with an effective vaccine to put an end to the pandemic. At the same time, researchers dedicated themselves to study, model and predict the evolution of the pandemic as well as investigate non-pharmaceutical intervention methods to control the spread of the virus \cite{Peng2020.02.16.20023465,early-transmission,zhao2020preliminary,wu2020nowcasting,Liu2020.01.25.919787,Shen2020.01.23.916726,kucharski2020early,data-driven,godio2020seir,GARBA2020108441,MANEVSKI2020108466,NARAIGH2020108496,NGONGHALA2020108364}.

On November 9th 2020 Pfizer and BioNTech \cite{pfizer,pfizer-publ} announced a vaccine candidate that successfully completed the clinical trials and is 90\% effective in preventing infection from the virus. Shortly after, two more vaccines were announced, namely the one by Moderna \cite{moderna} and the one by AstraZeneca \cite{astra}. All three vaccines got approved by the European Union \cite{ema-pfizer,ema-moderna,ema-astra} with a fourth vaccine getting granted a conditional marketing authorisation in March 2021, namely the Jcovden one (previously known as Janssen) \cite{ema-johnson}. Worldwide, certain states have made different decisions regarding available vaccines, as there are more vaccines available. Since the start of 2021, a vaccination roll-out commenced in most countries, including Ireland, taking into account the number of available vaccines and the level of risk different groups of people are considered to be in. As of 2022, first world countries have for the most part completed the first round of vaccinations and are in the process of administering booster shots. On the other hand, certain countries are expected to gain access to sufficient vaccines as late as 2023. \cite{padma2021covid}. Given this situation, a need that naturally arises is that of a way to determine the optimal vaccination strategy, especially given that the resources to reduce the severity of the pandemic are still limited \cite{ACUNAZEGARRA2021108614,RAO2021108621}. This idea was the main motivation for our study. 

Modelling, predicting and controlling the behaviour of epidemics has been a widely studied area \cite{booklenhart2007optimal,bailey1975mathematical}. A very prominent example is the pandemic influenza, a virus that caused an outbreak of severe pneumonia in 2009, commonly known as "swine-flu" \cite{chowell2009severe}. A number of models were developed to evaluate the implementation of mitigation strategies \cite{bowman2011evaluation,patel2005finding,ferguson-nature,Lee2012}, with a great focus on optimising these strategies. 

In this work, we introduce a new model that describes the transmission dynamics of the virus among different groups of the population in a more complete way compared to the well known SEIR model \cite{Peng2020.02.16.20023465,bacaer2011short,Mkendrick1925,anderson1979population} which is commonly used in the study of epidemics. The novelty of our approach lies in the introduction of new compartments to the model accounting for stages of the vaccination process as well as two different age groups. Using this model as our baseline, we explore a way to determine the optimal vaccination strategy, applying optimal control theory methods \cite{Fleming1975,ogata2010modern,kirk2004optimal,BOLTYANSKI1960464}. Using these tools, we studied the evolution of the pandemic in the Republic of Ireland, starting from January 2021, using data provided by \textit{Ireland's COVID-19 Data-Hub} \cite{datahub}. We obtained the optimal vaccination strategy based on estimates of the initial conditions, which is similar to the course of action taken by the state. However, we find that it is beneficial to start vaccinating people under the age of 65 in parallel with people in older age groups but in smaller numbers, as opposed to exclusively vaccinating the older population first. We compare our obtained strategy with the course of action taken by the state and comment on the benefits and disadvantages of both approaches.

In Section \ref{section_model} we introduce our model and express it as a system of ordinary differential equations. In Section \ref{section optimisation} we use optimal control theory techniques to obtain the optimal vaccination strategy, given our model. In Section \ref{section simulations} we produce two sets of simulations to compare the evolution of the pandemic with and without the vaccination strategy in place, under strict and loose restrictions. In Section \ref{section comparison} we present a comparative study of the strategy obtained through our optimisation against the strategy chosen by the government. Finally, in Section \ref{section discussion} we discuss our findings, the possible drawbacks of our method and extensions to our approach that are worth exploring. 

\section{Model}\label{section_model}

The model most commonly used in the study of epidemics was first introduced by \textit{M'Kendrick} in \cite{Mkendrick1925,anderson1979population} and is known as the SIR model. This model suggests that every member of a population can be considered to belong in one of the three compartments: Susceptible (S), Infectious (I) or Recovered (R). In the study of COVID-19, an extra compartment is usually added to the SIR model, namely the Exposed (E) compartment, consisting of the people who have been in contact with the virus yet have not developed any symptoms. This is known as the SEIR model \cite{Peng2020.02.16.20023465}. All compartments can be indexed by time (t), as they are expressing the number of individuals in each compartment for each moment in time.  In a closed population with no births and no deaths, the model can be expressed as a system of ordinary differential equations as follows:

\begin{eqnarray*}
	\frac{dS_t}{dt}&=& - \frac{\beta I_t}{N} S_t \\
	\frac{dE_t}{dt}&=& \frac{\beta I_t}{N} S_t -\sigma E_t \\ 
	\frac{dI_t}{dt}&=& \sigma E_t - \gamma I_t\\
	\frac{dR_t}{dt}&=& \gamma I_t \\
\end{eqnarray*}

where $N=S_t+E_t+I_t+R_t$ and it is constant since we are dealing with a closed population. $\beta$ is the infectious rate, which expresses the probability that a susceptible person gets infected by an infectious person, $\sigma$ is the incubation rate, meaning the rate at which an exposed individual becomes infectious and $\gamma$ is the recovery rate calculated as the inverse of the average time of infection. 

In this study we use a model that takes the idea of the SEIR model one step further, with additional compartments related to whether an individual has been vaccinated and whether the vaccination was effective, hence protecting the individual. The new set of states consists of the following: Susceptible not yet vaccinated (S); Received the vaccine; waiting for it to take effect (V); Received vaccine but it was not effective (N); Susceptible, refusing or unable to receive the vaccine (U); Exposed, infectious but still asymptomatic (E); Infectious symptomatic (I); Recovered or deceased (R); Protected from vaccine (P). 

\begin{figure}[!ht]
	\includegraphics[clip, trim=0.2cm 3.7cm 0.2cm 0.8cm,width=\linewidth]{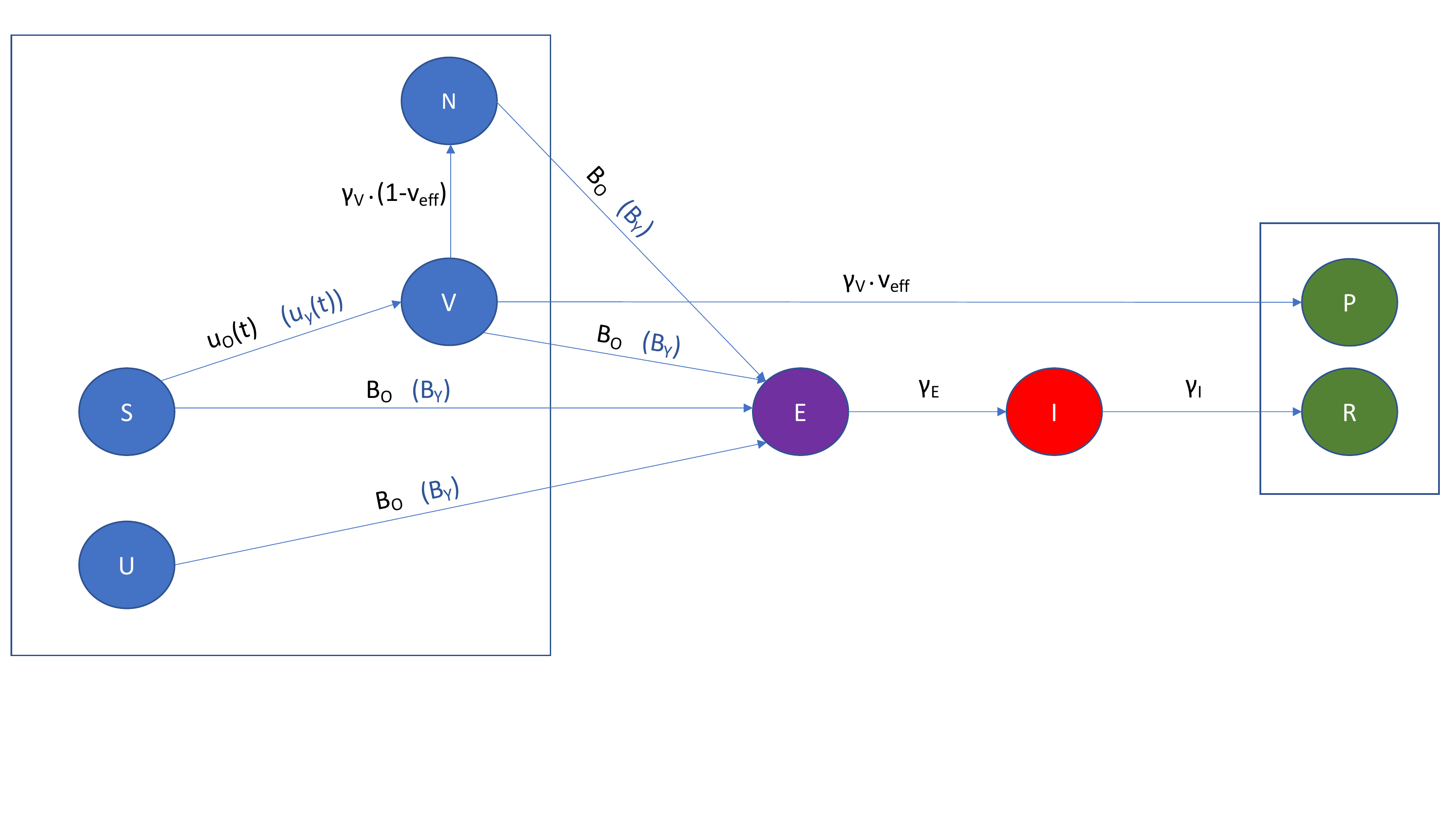}
	\caption[Vaccination Compartmental model]%
	{Vaccination Compartmental model \par \small The model compartments (states) are indicated by the coloured circles. Blue indicates states at which one might get infected by an exposed or infectious person and green indicates the states at which a person is considered to be safe from infection. The parameters describing the transmission rates are indicated on the arcs between the states. The models for the over and under 65 populations (indicated by O and Y subscripts) are identical with the only differences being the age specific control functions $u_O(t), u_Y(t)$ which express the percentage of the susceptible population to be vaccinated at each time point, and the rates of infection from an exposed or infectious individual $B_O =\frac{\beta_{OO} ( E_O(t)+I_O(t))}{T_O } +  \frac{\beta_{YO} ( E_Y(t)+I_Y(t) )}{T_Y } $ and $B_Y= \frac{\beta_{OY} ( E_O(t)+I_O(t) )}{T_O } +  \frac{\beta_{YY} ( E_Y(t)+I_Y(t) )}{T_Y }$.}
	\label{model_figure}
\end{figure}

A full table of the notations used throughout the paper can be found in \ref{appendix_notations}. In addition, we assume that a susceptible individual may get infected by exposed as well as infectious individuals. Furthermore, we consider two different age groups, those over 65 years old and those younger than 65. The reasons for this distinction can be found in \cite{mueller2020does}, the main ones being that 80\% of the hospitalised individuals are over the age of 65 and this age group has 23 times higher risk of death compared to those under 65. This results in two identical models, one for the younger population and one for the older which allows for the introduction of different parameters to account for the way the virus affects people of different ages. \figurename{ \ref{model_figure}} is a depiction of our model for an individual age cohort (for brevity we do not show the full model configuration). The state configuration is exactly the same in both cases with the difference being in the transition rates from one state to the other. 

We distinguish the states of the two models by adding an O or a Y as an subscript to the state name to indicate the older or younger than 65 age group respectively e.g. $S_O$, $V_Y$ etc. Since the states are time dependent, we can write them in function notation e.g. $S_O(t)$, $V_Y(t)$, however, for brevity, we will omit the $(t)$ in the equations that follow. We also use the notation $T_O, T_Y$ to indicate the total number of people in each age group which is assumed to not be time dependent. The two populations influence each other in the sense that a susceptible person in any of the two groups may be infected by an exposed or infectious person from either group. For brevity, we will indicate the older than 65 age group as o65 and the younger than 65 age group as y65. We will be using upper case Roman characters to indicate the model states and lower case Greek characters to indicate model parameters.

The model is expressed as a system of ordinary differential equations (ODEs), each of them describing the evolution through time of one of the states. The ODEs describing the dynamics of our model are given by:

\allowdisplaybreaks
\begin{eqnarray}\label{stateequations}
\frac{dS_O}{dt} &=& -\left(  \frac{\beta_{OO} ( E_O+ I_O )}{T_O } +  \frac{\beta_{YO}( E_Y+I_Y)}{T_Y } \right)   S_O  -u_{O} S_O\nonumber\\
\frac{dS_Y}{dt} &=& -\left(  \frac{\beta_{OY} ( E_O+I_O)}{T_O } +  \frac{\beta_{YY} ( E_Y+I_Y )}{T_Y } \right)   S_Y  -u_{Y}  S_Y \nonumber\\
\frac{dV_O}{dt} &=& u_{O} S_O - \left(  \frac{\beta_{OO} ( E_O+I_O )}{T_O } +  \frac{\beta_{YO} ( E_Y+I_Y )}{T_Y } \right)V_O  \nonumber \\ & &  - (1-\alpha_V)  \gamma_V  V_O -\gamma_V \alpha_V  V_O\nonumber\\
\frac{dV_Y}{dt} &=& u_{Y}  S_Y - \left(  \frac{\beta_{OY} ( E_O+I_O )}{T_O } +  \frac{\beta_{YY} ( E_Y+I_Y )}{T_Y } \right)  V_Y \nonumber \\ & &  - (1-\alpha_V)  \gamma_V  V_Y -\gamma_V  \alpha_V   V_Y\nonumber\\
\frac{dN_O}{dt} &=&(1-\alpha_V)\gamma_V V_O - \left(  \frac{\beta_{OO} ( E_O+I_O )}{T_O } +  \frac{\beta_{YO} ( E_Y+I_Y )}{T_Y } \right) N_O \nonumber\\
\frac{dN_Y}{dt} &=&(1-\alpha_V)  \gamma_V  V_Y - \left(  \frac{\beta_{OY} ( E_O+I_O)}{T_O } +  \frac{\beta_{YY} ( E_Y+I_Y )}{T_Y } \right)  N_Y\nonumber \\
\frac{dU_O}{dt} &=&  -\left(  \frac{\beta_{OO} ( E_O+I_O)}{T_O } +  \frac{\beta_{YO}( E_Y+I_Y )}{T_Y } \right)   U_O\nonumber\\
\frac{dU_Y}{dt} &=&   -\left(  \frac{\beta_{OY}\cdot ( E_O+I_O)}{T_O } +  \frac{\beta_{YY}\cdot ( E_Y+I_Y )}{T_Y } \right)   U_Y \nonumber\\
\frac{dE_O}{dt} &=& \left(  \frac{\beta_{OO}   ( E_O+I_O )}{T_O } +  \frac{\beta_{YO} ( E_Y+I_Y )}{T_Y } \right)  \left( S_O +V_O  + N_O    + U_O\right)  \nonumber\\ & & 
- \gamma_E E_O \nonumber\\
\frac{dE_Y}{dt} &=& \left(  \frac{\beta_{OY}( E_O+I_O )}{T_O} +  \frac{\beta_{YY} ( E_Y+I_Y )}{T_Y } \right)   \left( S_Y+V_Y  + N_Y + U_Y\right)   \nonumber\\ & &
- \gamma_E E_Y\nonumber\\
\frac{dI_O}{dt} &=& \gamma_E E_O -\gamma_I I_O \nonumber\\
\frac{dI_Y}{dt} &=& \gamma_E E_Y -\gamma_I I_Y \nonumber\\
\frac{dR_O}{dt} &=&\gamma_I I_O  \nonumber\\
\frac{dR_Y}{dt} &=&\gamma_I I_Y  \nonumber\\
\frac{dP_O}{dt} &=&\gamma_V \alpha_V \cdot  V_O\nonumber\\
\frac{dP_Y}{dt} &=&\gamma_V \alpha_V \cdot  V_Y
\end{eqnarray}

where $T_O,T_Y$ denote the total number of people in the over and under 65 years old populations respectively. These populations are divided in those willing and those unwilling to get vaccinated ($S_O,S_Y,U_O,U_Y$). This division happens based on estimates of the percentages of people who refuse the vaccine in each age group. The only option for individuals in the $U$ compartments, apart from staying there if the pandemic ends early, is to eventually get infected and be moved to the exposed compartments $E$. 

The terms $u_O, u_Y$ are the control functions and represent the percentage of the old and young population respectively to be vaccinated at each time point, which means that they are also time dependent. Each control function is applied to the respective susceptible population ($S_O, S_Y$), moving that percentage of the population to the respective vaccinated compartment ($V_O, V_Y$). The optimal values for $u_O, u_Y$ are obtained through the application of the Optimal Control techniques \cite{Fleming1975,neustadt1976optimization} discussed in the next section. Any person in the states $(S,V,U,N)$ is considered susceptible to getting infected by an infectious or exposed person $(E,I)$. That is because they are either not yet chosen for vaccination, they received the vaccine and are waiting for it to take effect, they received the vaccine and it was ineffective, or they chose not to get vaccinated. The terms $\beta_{ij}, i,j \in \{O,Y \} $ describe the transmission rates of the infection between the age groups and they are analogous to the $\beta$ parameter in the standard SEIR model. Specifically: 
\begin{itemize}
	\item $\beta_{OO} = $ rate at which an o65 susceptible person becomes infected by an o65 exposed or o65 infected person
	\item $\beta_{YO} = $ rate at which an o65 susceptible person becomes infected by an y65 exposed or y65 infected person
	\item $\beta_{OY} = $ rate at which a y65 susceptible person becomes infected by an o65 exposed or o65 infected person
	\item $\beta_{YY} = $ rate at which a y65 susceptible person becomes infected by an y65 exposed or y65 infected person
\end{itemize}
There are three rates taken into consideration in the model, the first being $\gamma_E$ which is the rate at which an exposed person becomes symptomatic. This rate can be calculated as the inverse of the mean holding time to develop symptoms and become infectious. The next rate is $\gamma_I$ which is the rate at which an infectious person recovers and can be calculated as the inverse of the mean time to recovery. The final rate considered in the model is $\gamma_V$ which is the rate at which vaccination becomes effective and can be calculated as the inverse of the mean holding time until protected from the vaccine. The last parameter influencing the model is $\alpha_V$ which is the vaccine effectiveness. This is expressed as the percentage of people who become protected from the vaccine, out of those who have received it. 

\section{Optimal control} \label{section optimisation}
Optimal control theory \cite{Fleming1975,neustadt1976optimization,kirk2004optimal} is the study of strategies to obtain the control function that optimises a certain objective over a time horizon, possibly under suitable constraints. These types of techniques have been widely adopted to biological systems in general \cite{booklenhart2007optimal} but also more specifically to obtain optimal vaccination strategies when dealing with viruses and epidemics \cite{optimal-epidemics,Lee2012,sharomi2017optimal,kirschner1997optimal}. More specifically, \textit{Pontryagin's Maximum Principle} \cite{pontryagin2018mathematical,BOLTYANSKI1960464,lee1967foundations} is the main tool that is used when dealing with a problem whose dynamics are described by a set of Ordinary Differential Equations.

An optimisation problem can be expressed as the problem of minimising an objective (cost) functional under certain constraints. Let $f(t,x,u)$ denote the objective functional. Also, let $g(t,x,u)$ denote the state equation (or set of state equations) of our system. Using $f$ and $g$ we can form a Hamiltonian function \cite{kirk2004optimal,ogata2010modern} as follows:
\begin{equation*}
H(t,x,u,\lambda ) =  f(t,x,u)+ \lambda(t) g(t,x,u)
\end{equation*}

where $\lambda$ is a continuous function of time$(t)$. For simplicity, we will write $\lambda(t)$ as just $\lambda$ in the expressions that follow. 

For the COVID-19 control problem and our spcific model, we can express the goal of our optimisation as the minimisation of the number of infectious individuals, at a minimal cost via vaccination, within a certain time frame $[0,T]$. That goal can be expressed with the help the following objective functional to be minimised:

\begin{equation}
\mathcal{F}\left( U(t)\right) = \int_{0}^{T} \left[ I_O(t)+ I_Y(t) + \frac{W_O}{2}u_O^2(t) + \frac{W_Y}{2}u_Y^2(t) \right] dt 
\end{equation}

where $U(t)=\left( u_O(t),u_Y(t)\right) $ and $W_O, W_Y \geq 0$ are the age specific weight constants enforcing the severity of the optimisation constraints.  This type of choice of the objective functional as a non-negative increasing function of the control inputs, is common practice in optimal control applications \cite{Lee2012,sharomi2017optimal,kirk2004optimal}. The control functions are squared in order to ensure the convexity of the functional, guaranteeing a sufficient condition for the existence of an optimal solution (to be proved later in detail). Consider $ X(t)= $ $\left( \right.S_O(t),$ $V_O(t),$ $N_O(t) ,$ $U_O(t),$  $E_O(t),$ $I_O(t),$ $R_O(t),$  $P_O(t),$ $S_Y(t),$ $V_Y(t),$ $N_Y(t),$  $U_Y(t),$ $E_Y(t),$ $I_Y(t),$ $R_Y(t) ,$ $P_Y(t) \left.\right) $. We are looking for the optimal pair of solutions $\left( U^*(t),X^*(t)\right) $, i.e. the optimal control $U^*$ and the corresponding trajectory $X^*$ when $U^*$ is applied, such that

\begin{equation} \label{minimization of functional}
\mathcal{F}(U^*(t))=\min_\Omega \mathcal{F}(U(t))
\end{equation}

where $\Omega = \left\lbrace U(t) \in L^2(O,T)^2 \Vert a \leq u_O(t),u_Y(t)\leq b, t \in [0,T]\right\rbrace $, $a$ and $b$ are the upper and lower bounds for the control function and can usually be expressed as real values, and $T$ is the time horizon for our optimisation. 

There are a few different approaches that can be taken in defining the constrained optimisation problem and are worth mentioning. We can include weighting factors to the infectious compartments ($I_O,I_Y$) in the objective functional $\mathcal{F}$ that correspond to the mortality rates for each age group or more factors that generally model the cost that large number of infections can have in the healthcare system. Additionally, when it comes to bounding the control functions, we can chose a constraint of the form $a \leq u_O(t)+u_Y(t)\leq b$ which can be interpreted as the total percentage of vaccinations being bounded as opposed to bounding them per age group. This kind of constraint would result in an extra term in the objective functional and the Hamiltonian function having the form $\lambda \cdot (u_O+u_Y)$.

The following theorem describes the main conditions that when satisfied can lead us to the optimal control that solves the minimisation problem.

\begin{theorem}[Pontryagin's maximum principle] \label{pontryagin}
	Given the Hamiltonian 
	\begin{equation*}
	H(t,x,u,\lambda ) =  f(t,x,u)+ \lambda g(t,x,u),
	\end{equation*}
	then the following conditions are satisfied by the optimal control $u^*\in \Omega= \left\lbrace U(t) \in L^2(O,T)^2 \Vert a \leq u_O(t),u_Y(t)\leq b, t \in [0,T]\right\rbrace$:
	\begin{eqnarray*}
		\frac{\partial H}{\partial u} \left\lbrace \begin{matrix} =&0 & \text{ at } u^* \text{ if } &a<u^*<b \\ \geq & 0 & \text{ at } u^* \text{ if }& u^*=a\\ \leq & 0 & \text{ at } u^* \text{ if }& u^*=b\end{matrix} \right. & & \text{Optimality Condition}\\
		\dot{\lambda } = - \frac{\partial H}{\partial x}  & &\text{Adjoint Equation} \\
		\lambda(T) =0 & &\text{Transversality Condition} \\
		\dot{x} =g(t,x,u), & x(0)=x_0   &\text{Dynamics of state equation}
	\end{eqnarray*}
	where $\dot{\lambda}=\partial \lambda / \partial t$ and  $\dot{x}=\partial x / \partial t$.
\end{theorem}

In order to apply the maximum principle we first need to define the Hamiltonian function, omitting for brevity the time dependence from the state notations:

\begin{eqnarray*}\label{hamiltonian}
	H& =&  \left[ I_O + I_Y + \frac{W_O}{2}u_O^2 + \frac{W_Y}{2}u_Y^2 \right]  \\
	& + & \lambda_{S_O} \left\lbrace -\left(  \frac{\beta_{OO} ( E_O+I_O )}{T_O } +  \frac{\beta_{YO} ( E_Y+I_Y )}{T_Y } \right)   S_O -u_{O}  S_O \right\rbrace \\
	& + & \lambda_{S_Y} \left\lbrace -\left(  \frac{\beta_{OY}( E_O+I_O )}{T_O } +  \frac{\beta_{YY}( E_Y+I_Y)}{T_Y } \right)   S_Y  -u_{Y} S_Y \right\rbrace \\
	&+& \lambda_{V_O} \left\lbrace u_{O}  S_O - \left(  \frac{\beta_{OO} ( E_O+I_O )}{T_O } +  \frac{\beta_{YO} ( E_Y+I_Y )}{T_Y } \right)   V_O  - \gamma_V  V_O \right. \\ & &\left.-\alpha_V   V_O\right\rbrace \\
	&+& \lambda_{V_Y}\left\lbrace u_{Y}  S_Y - \left(  \frac{\beta_{OY} ( E_O+I_O)}{T_O } +  \frac{\beta_{YY}( E_Y+I_Y )}{T_Y } \right)   V_Y  - \gamma_V V_Y\right. \\ & &\left. -\alpha_V  V_Y \right\rbrace \\
	&+& \lambda_{N_O} \left\lbrace \gamma_V  V_O - \left(  \frac{\beta_{OO}( E_O+I_O)}{T_O } +  \frac{\beta_{YO} ( E_Y+I_Y )}{T_Y } \right)   N_O \right\rbrace \\
	&+& \lambda_{N_Y} \left\lbrace \gamma_V V_Y - \left(  \frac{\beta_{OY}( E_O+I_O )}{T_O } +  \frac{\beta_{YY} ( E_Y+I_Y )}{T_Y } \right)   N_Y \right\rbrace \\
	&-& \lambda_{U_O} \left(  \frac{\beta_{OO} ( E_O+I_O)}{T_O } +  \frac{\beta_{YO} ( E_Y+I_Y)}{T_Y} \right) U_O\\
	&-& \lambda_{U_Y} \left(  \frac{\beta_{OY} ( E_O+I_O)}{T_O } +  \frac{\beta_{YY}( E_Y+I_Y)}{T_Y } \right)   U_Y\\
	&+& \lambda_{E_O} \left\lbrace  \left(  \frac{\beta_{OO} ( E_O+I_O)}{T_O } +  \frac{\beta_{YO} ( E_Y+I_Y )}{T_Y } \right)  \left( S_O +V_O + N_O + U_O\right)  \right. \\ & & \left. 
	- \gamma_E E_O \right\rbrace \\
	&+& \lambda_{E_Y} \left\lbrace  \left(  \frac{\beta_{OY}( E_O+I_O )}{T_O } +  \frac{\beta_{YY} ( E_Y+I_Y)}{T_Y } \right)   \left( S_Y +V_Y + N_Y + U_Y\right) \right. \\ & & \left.
	- \gamma_E E_Y \right\rbrace \\
	&+& \lambda_{I_O} \left\lbrace  \gamma_E E_O -\gamma_I I_O\right\rbrace  \\
	&+& \lambda_{I_Y} \left\lbrace   \gamma_E E_Y -\gamma_I I_Y \right\rbrace \\
\end{eqnarray*}
where $\Lambda:=$ $\left\lbrace \right.\lambda_{S_O},$ $\lambda_{S_Y},$ $\lambda_{V_O}, $ $\lambda_{V_Y},$ $\lambda_{N_O},$ $\lambda_{N_Y},$ $\lambda_{U_O},$ $\lambda_{U_Y},$ $\lambda_{E_O},  $
$ \lambda_{E_Y},$ $\lambda_{I_O},$ $\lambda_{I_Y} \left.\right\rbrace$ is a set of continuous functions of time ($\Lambda(t)$), given the states and controls. These functions have a key role in our optimisation technique.

Applying each of the conditions of the maximum principle (\ref{pontryagin}) we obtain the following:  
\begin{itemize}
	\item From the optimality condition, the optimal control strategies are given by:
	\begin{equation}\label{optimal controlO2}
	u_O^*(t)= \min \left\lbrace \max \left \lbrace a,\frac{S_O(t)}{W_O} \left( \lambda_{S_O}(t) -\lambda_{V_O} (t)\right) \right \rbrace, b \right\rbrace 
	\end{equation}
	
	\begin{equation}\label{optimal controlY2}
	u_Y^*(t)= \min \left\lbrace \max \left \lbrace a,\frac{S_Y(t)}{W_Y}\left( \lambda_{S_Y}(t) -\lambda_{V_Y} (t)\right)  \right \rbrace, b \right\rbrace 
	\end{equation}
	
	\item Adjoint Equation: $\dot{\Lambda } = - \partial H / \partial X$
	which results in the following system of ODEs:
\end{itemize}

\begin{eqnarray*}\label{adjointequations}
	\dot{\lambda}_{S_O} &=& \left( \frac{\beta_{OO} ( E_O+I_O )}{T_O } +  \frac{\beta_{YO}( E_Y+I_Y )}{T_Y }    +u_{O}\right)  \lambda_{S_O} -u_O \lambda_{V_O}\\
	&-&  \left(  \frac{\beta_{OO} ( E_O+I_O)}{T_O } +  \frac{\beta_{YO} ( E_Y+I_Y )}{T_Y } \right)  \lambda_{E_O}\\
	\dot{\lambda}_{S_Y} &=& \left( \frac{\beta_{OY} ( E_O+I_O)}{T_O } +  \frac{\beta_{YY} ( E_Y+I_Y)}{T_Y }    +u_{Y} \right) \lambda_{S_Y}- u_Y\lambda_{V_Y}\\
	& - & \left(  \frac{\beta_{OY}( E_O+I_O )}{T_O } +  \frac{\beta_{YY} ( E_Y+I_Y )}{T_Y } \right)  \lambda_{E_Y}\\
	\dot{\lambda }_{V_O} &=& \left( \frac{\beta_{OO}( E_O+I_O)}{T_O } +  \frac{\beta_{YO} ( E_Y+I_Y )}{T_Y }  + \gamma_V\right) \lambda_{V_O} - (1-\alpha_V)  \gamma_V \lambda_{N_O}\\ 
	& -& \left(  \frac{\beta_{OO}( E_O+I_O )}{T_O } +  \frac{\beta_{YO} ( E_Y+I_Y)}{T_Y } \right) \lambda_{E_O}\\
	\dot{\lambda }_{V_Y} &=& \left( \frac{\beta_{OY}( E_O+I_O )}{T_O } +  \frac{\beta_{YY} ( E_Y+I_Y )}{T_Y }  +\gamma_V \right) \lambda_{V_Y} -(1-\alpha_V)   \gamma_V \lambda_{N_Y}\\ 
	& -&\left(  \frac{\beta_{OY} ( E_O+I_O )}{T_O } +  \frac{\beta_{YY} ( E_Y+I_Y )}{T_Y } \right)  \lambda_{E_Y}\\
	\dot{\lambda}_{N_O} &=& \left(  \frac{\beta_{OO}( E_O+I_O )}{T_O } +  \frac{\beta_{YO}( E_Y+I_Y)}{T_Y } \right) (\lambda_{N_O}-\lambda_{E_O}) \\
	\dot{\lambda}_{N_Y} &=& \left(  \frac{\beta_{OY} ( E_O+I_O )}{T_O } +  \frac{\beta_{YY}( E_Y+I_Y)}{T_Y } \right)  (\lambda_{N_Y}-\lambda_{E_Y}) \\ 
	\dot{\lambda}_{U_O} &=&  \left(  \frac{\beta_{OO}( E_O+I_O )}{T_O } +  \frac{\beta_{YO}(E_Y+I_Y)}{T_Y} \right)  (\lambda_{U_O}-\lambda_{E_O})\\
	\dot{\lambda}_{U_Y} &=&  \left(  \frac{\beta_{OY}(E_O+I_O)}{T_O } +  \frac{\beta_{YY}(E_Y+I_Y)}{T_Y } \right) (\lambda_{U_Y}-\lambda_{E_Y})\\
	\dot{\lambda}_{E_O} &=& \frac{\beta_{OO} S_O }{T_O}\left( \lambda_{S_O}-\lambda_{E_O}\right) + \frac{\beta_{OO} V_O }{T_O} \left( \lambda_{V_O}-\lambda_{E_O}\right) \\
	& + & \frac{\beta_{OO} N_O }{T_O} \left( \lambda_{N_O}-\lambda_{E_O}\right) + \frac{\beta_{OO} U_O }{T_O} \left( \lambda_{U_O}-\lambda_{E_O}\right) \\
	&+&   \frac{\beta_{OY} S_Y }{T_O}\left( \lambda_{S_Y}-\lambda_{E_Y}\right) 
	+ \frac{\beta_{OY} YV_t }{T_O} \left( \lambda_{V_Y}-\lambda_{E_Y}\right) \\&+& \frac{\beta_{OY} N_Y }{T_O} \left( \lambda_{N_Y}-\lambda_{E_Y}\right) 
	+ \frac{\beta_{OY} U_Y }{T_O} \left( \lambda_{U_Y}-\lambda_{E_Y}\right) \\
	&+& \gamma_E \lambda_{E_O}  - \gamma_E \lambda_{I_O}\\
	\dot{\lambda}_{E_Y} &=& \frac{\beta_{YO}S_O }{T_Y}\left( \lambda_{S_O}-\lambda_{E_O}\right) + \frac{\beta_{YO} V_O }{T_Y} \left( \lambda_{V_O}-\lambda_{E_O}\right) \\
	& + & \frac{\beta_{YO}N_O }{T_Y} \left( \lambda_{N_O}-\lambda_{E_O}\right)+ \frac{\beta_{YO} U_O }{T_Y} \left( \lambda_{U_O}-\lambda_{E_O}\right) \\
	&+& \frac{\beta_{YY} S_Y }{T_Y}\left( \lambda_{S_Y}-\lambda_{E_Y}\right) 
	+ \frac{\beta_{YY}V_Y }{T_Y} \left( \lambda_{V_Y}-\lambda_{E_Y}\right)\\ &+&\frac{\beta_{YY} N_Y }{T_Y} \left( \lambda_{N_Y}-\lambda_{E_Y}\right) 
	+ \frac{\beta_{YY} U_Y }{T_Y} \left( \lambda_{U_Y}-\lambda_{E_Y}\right) \\ &+& \gamma_E  \lambda_{E_Y}  - \gamma_E  \lambda_{I_Y} \\
	\dot{\lambda}_{I_O} &=& \frac{\beta_{OO}S_O }{T_O}\left( \lambda_{S_O}-\lambda_{E_O}\right) + \frac{\beta_{OO}V_O }{T_O} \left( \lambda_{V_O}-\lambda_{E_O}\right) \\
	& + & \frac{\beta_{OO} N_O }{T_O} \left( \lambda_{N_O}-\lambda_{E_O}\right)+ \frac{\beta_{OO} U_O }{T_O} \left( \lambda_{U_O}-\lambda_{E_O}\right) \\
	&+&  \frac{\beta_{OY} S_Y }{T_Y}\left( \lambda_{S_Y}-\lambda_{E_Y}\right) 
	+\frac{\beta_{OY} V_Y }{T_O} \left( \lambda_{V_Y}-\lambda_{E_Y}\right)\\ &+&\frac{\beta_{OY}N_Y }{T_O} \left( \lambda_{N_Y}-\lambda_{E_Y}\right)
	+ \frac{\beta_{OY} U_Y }{T_O} \left( \lambda_{U_Y}-\lambda_{E_Y}\right) \\ &+& \gamma_I  \lambda_{I_O} -1 \\
	\dot{\lambda}_{I_Y} &=& \frac{\beta_{YO} S_O }{T_Y}\left( \lambda_{S_O}-\lambda_{E_O}\right) + \frac{\beta_{YO}V_O }{T_Y} \left( \lambda_{V_O}-\lambda_{E_O}\right) \\
	& + & \frac{\beta_{YO} N_O }{T_Y} \left( \lambda_{N_O}-\lambda_{E_O}\right)+ \frac{\beta_{YO} U_O }{T_Y} \left( \lambda_{U_O}-\lambda_{E_O}\right) \\
	&+& \frac{\beta_{YY} S_Y }{T_Y}\left( \lambda_{S_Y}-\lambda_{E_Y}\right) 
	+ \frac{\beta_{YY} V_Y }{T_Y} \left( \lambda_{V_Y}-\lambda_{E_Y}\right)\\ &+&\frac{\beta_{YY}N_Y }{T_Y} \left( \lambda_{N_Y}-\lambda_{E_Y}\right) 
	+ \frac{\beta_{YY} U_Y }{T_Y} \left( \lambda_{U_Y}-\lambda_{E_Y}\right) \\ &+& \gamma_I  \lambda_{I_Y} -1
\end{eqnarray*}
where $\dot{\lambda}$ denotes $\partial \lambda / \partial t$

\begin{itemize}
	\item The transversality condition: $\Lambda(T)=0$ gives us a terminal condition for the system of adjoint ODEs. This condition implies that the adjoint system should be solved backwards \cite{butcher2008numerical} as opposed to the state system. 
	\item The dynamics of the system are described by the system of state equations (\ref{stateequations}).
\end{itemize}

To ensure the existence of a solution to our optimisation problem, we need to examine the sufficient conditions as stated in \cite{Fleming1975}, as also used in \cite{gaff2009optimal}. The relevant conditions and the corresponding existence proof for our model can be found in Appendix \ref{appendix_conditions}.

Having obtained all the necessary optimality conditions and resulting equations,  Algorithm \ref{algorithm} describes the steps that lead to the computation of an optimal control trajectory.  

 \begin{algorithm}
	\SetAlgoLined
	\KwResult{Optimal control and trajectory $U^*, X^*$}
	\KwIn{Initial state $x_0$, Parameter values, Time horizon, Initial control, Convergence threshold $\delta$ }
	Solve the system of state ODEs forwards in time to determine the state\; 
	Solve the system of adjoint ODEs backwards in time\;
	Update the controls using (\ref{optimal controlO2}) and (\ref{optimal controlY2}) \; 
	Compute the error term $ \frac{\Vert x_k - x_{k-1} \Vert_1}{\Vert x_{k-1} \Vert_1} $ \;
	If the error term is smaller than some predetermined threshold $\delta$ then extract the optimal control and optimal trajectory. Otherwise repeat the process until convergence.  
	\caption{ \label{algorithm} Determining the optimal control}
\end{algorithm}

\section{Simulation Setup and Results} \label{section simulations}
We used our model to estimate the evolution of the virus in the Republic of Ireland, starting from January 2021 using data provided by \cite{datahub} regarding the population of the country, the number of exposed, infected and recovered people in both age groups. Additionally, we applied the optimal control strategy on the same data to obtain the effect that the vaccination would have on the evolution of the pandemic, were it adopted by the state. 

Table \ref{initial states} includes the initial values of the states $E_O, I_O, R_O, E_Y, I_Y, R_Y$. The rest of the states can be calculated with the use of this information and some of the parameter values. The total susceptible people in each age group are the people not yet vaccinated along with the people refusing or unable to receive the vaccine ($\textbf{S}:=S+U$). This total of susceptibles for each age group can be calculated by: $\textbf{S}= T - (E+I+R)$. The separation between not yet vaccinated ($S$) and unwilling to get vaccinated ($U$) can be expressed as a percentage for each population, representing the refusal rate for the vaccine. Let $r_O$ and $r_Y$ be these percentages in the over 65 and under 65 populations respectively. Then $S_Y=(1-r_Y)\textbf{S}$, $U_Y=r_Y\textbf{S}$, $S_O=(1-r_O)\textbf{S}$, $U_O=r_O\textbf{S}$.

\begin{table}[ht]
	\centering
	\begin{tabular}{|l|r|} 
		\hline 
		State (Compartment) & Number of people \\
		\hline
		Total y65 ($T_Y$)& 4000000 \\
		Total o65 ($T_O$) & 900000 \\
		Y65 exposed ($E_Y$)& 2000 \\
		O65 exposed ($E_O$)& 200 \\
		Y65 infected ($I_Y$)& 2000 \\
		O65 infected ($I_O$)& 200 \\
		Y65 recovered ($R_Y$)& 200000 \\ 
		O65 recovered ($R_O$)& 100000 \\
		\hline
	\end{tabular}
	\caption{\label{initial states} Initial values for the model states }
\end{table} 

We produce two different simulations, to study the evolution of the pandemic, following the steps described in Algorithm \ref{algorithm}. In the first case we assume strict measures are in place, meaning that the transmission rates of the virus among individuals is low. In the second case we assume minimal restrictions are in place, hence the transmission rates are much higher. The parameters used in both simulations are listed in Table \ref{initial parameters}.

\begin{table}[ht]
	\centering
	\begin{tabular}{|l|c|c|c|} 
		\hline
		Parameter Description & Symbol& Case 1& Case 2\\
		\hline
		Mean no of o65 infected by an o65 &$R0_{OO}$ &1.2&8\\
		Mean no of y65 infected by a y65 &$R0_{YY}$ &1.2&8\\
		Mean no of y65 infected by an o65 &$R0_{OY}$  &0.9&4\\
		Mean no of o65 infected by a y65 &$R0_{YO}$ &0.9&3\\
		\hline
		Mean holding time exposed in days &$t_E$ &\multicolumn{2}{c|}{6.6}\\
		Mean holding time infected in days & $t_I$ &\multicolumn{2}{c|}{7.4}\\
		Mean holding time vaccinated in days & $t_V$ & \multicolumn{2}{c|}{14} \\
		Vaccine effectiveness &$\alpha_V$ & \multicolumn{2}{c|}{0.9}\\
		Upper bound for control function &$b$ & \multicolumn{2}{c|}{0.3}\\
		Percentage of o65 people refusing the vaccine & $r_O$ &\multicolumn{2}{c|}{0.07}\\
		Percentage of y65 people refusing the vaccine &$r_Y$ & \multicolumn{2}{c|}{0.21}\\
		Weight constant relating to o65 group &$W_O$ & \multicolumn{2}{c|}{$10^{11}$} \\
		Weight constant relating to y65 group &$W_Y$ & \multicolumn{2}{c|}{$10^{11}$}\\
		\hline
	\end{tabular} 
	\caption{\label{initial parameters}Parameter values for both simulations. The mean numbers of infections from each infected person within each age group and between age groups vary for the two simulations, while for the rest of the parameters we used the same values. Case 1 is representative of a situation where strict measures are in place and the results of the simulation can be found in \subsectionautorefname{ \ref{section tight restrictions}}. Case 2 is an example of a situation with minimal restrictions and is explored in \subsectionautorefname{ \ref{section loose restrictions}}.}
\end{table}

The three parameters $t_E,t_I,t_V$, describe the mean holding times in each of the states E, I and V. They are the inverses of the rates we have in our model (\ref{stateequations}), namely $\gamma_E, \gamma_I, \gamma_V$. Also, particularly interesting are the four $RO$ parameters, that are related to how quickly the infection is transmitted between the members of the entire population. More specifically, these parameters are directly related to the terms $\beta_{ij}$ in our model which describe the transmission rates, with the help of the three mean holding times:
\begin{eqnarray}\label{betar0}
\beta_{OO}  =  \frac{R0_{OO}}{t_E+t_I+t_V}, \beta_{YO} =\frac{R0_{YO}}{t_E+t_I+t_V}, \nonumber\\
\beta_{OY}  = \frac{R0_{OY}}{t_E+t_I+t_V}, \beta_{YY} =\frac{R0_{YY}}{t_E+t_I+t_V}
\end{eqnarray}
Two different sets of $R0$ numbers are used to produce our two simulations and give an insight to the evolution of the pandemic under different levels of restrictions. 

The values we use for the parameters are estimates specific to COVID-19 (mean holding times) \cite{cereda2020early}, or based on estimates used in similar studies of other infectious diseases such as (weight constants) \cite{Lee2012}. The upper bound for the control function ($b$) is chosen arbitrarily, however it is not a value that is ever obtained for either control function, thanks to the high weight constants used in the objective functional that model the cost of the vaccination and the severity of the constraints. An alternative approach would be for $b$ to bound the sum of the control functions as opposed to each of them individually. The weight constants $W_O, W_Y$ are chosen to be equal because we assume no difference in the cost of providing the vaccination to people in different age groups. However, the two groups can be weighed differently, explaining the cost of bringing the vaccination to remote locations, nursing homes etc. The refusal percentages for both age groups ($r_O, r_Y$) we used are based on a survey \cite{survey}, the results of which suggest that $77\% $ of the overall population of Ireland are willing to receive the vaccine and the the willingness is stronger in the over 65 population, with $93\%$ of the group intending to get vaccinated. We included the percentage of people unsure about the vaccination ($15\%$ overall) in the refusal percentage ($6\%$) along with the people who claimed to already be against receiving the vaccination. Finally, as we mentioned in \sectionautorefname{ \ref{section optimisation}}, it is possible for one more set of parameters to be included in our model and more specifically the optimisation function, weighing the number of infectious individuals by the mortality rate for each age group.

\subsection{\label{section tight restrictions}Evolution of pandemic under tight restrictions}

We simulate the evolution of the pandemic using the initial state values given in  \tableautorefname{ \ref{initial states}} and the parameter values given in \tableautorefname{ \ref{initial parameters}}. Specifically, the $R0$ numbers we use are given in the column named Case 1 and are an example of the way the virus is evolving under a strict restriction policy. In that case, the average number of infections within each age group is 1.2 and between age groups 0.9. 

\begin{figure}[ht]
	\begin{subfigure}{.45\textwidth}
		\includegraphics[width=\linewidth,height=4cm]{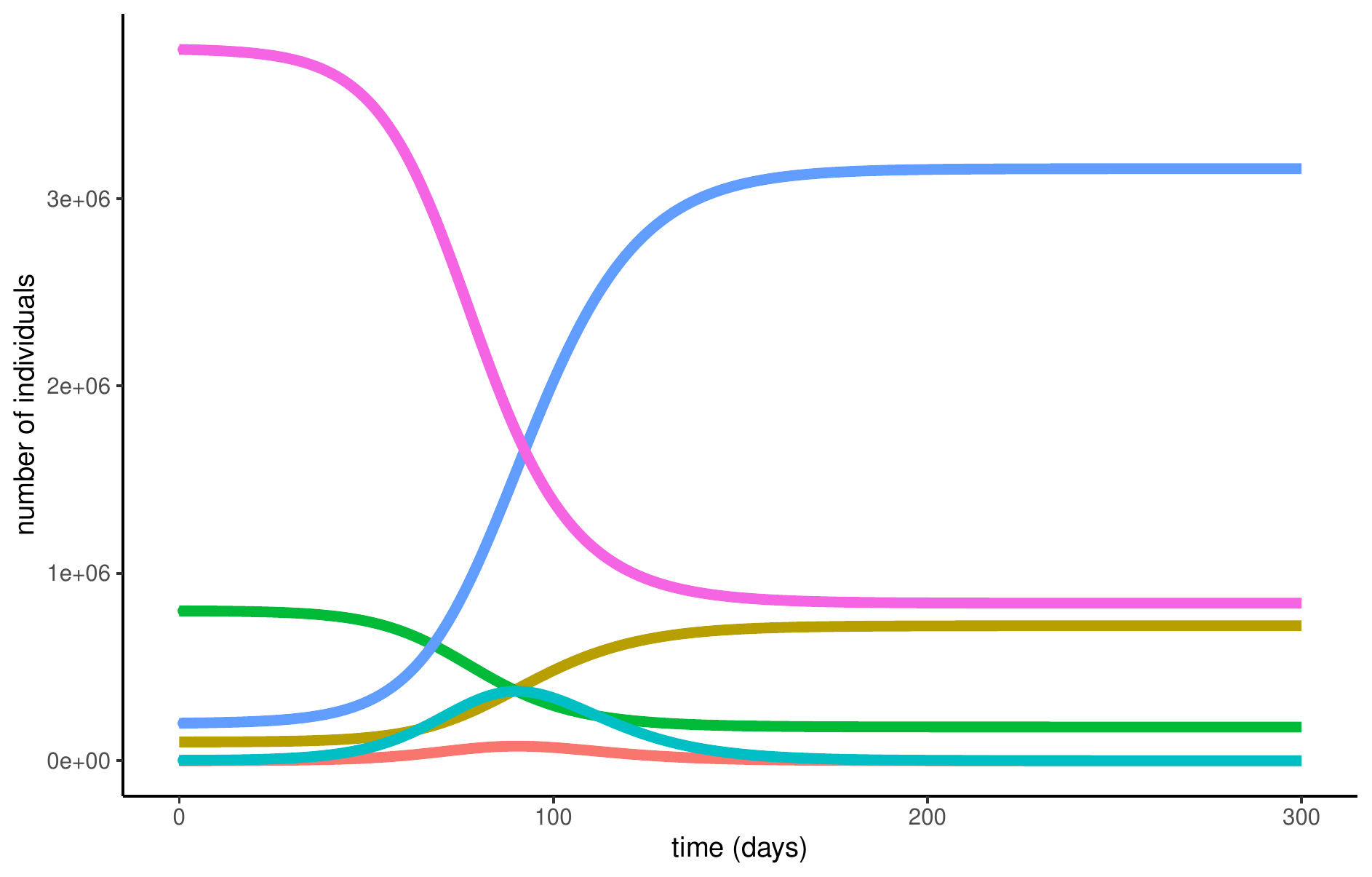}
		\caption{Evolution of states without the vaccine}
	\end{subfigure}%
	\begin{subfigure}{.55\textwidth}
		\includegraphics[width=\linewidth,height=4cm]{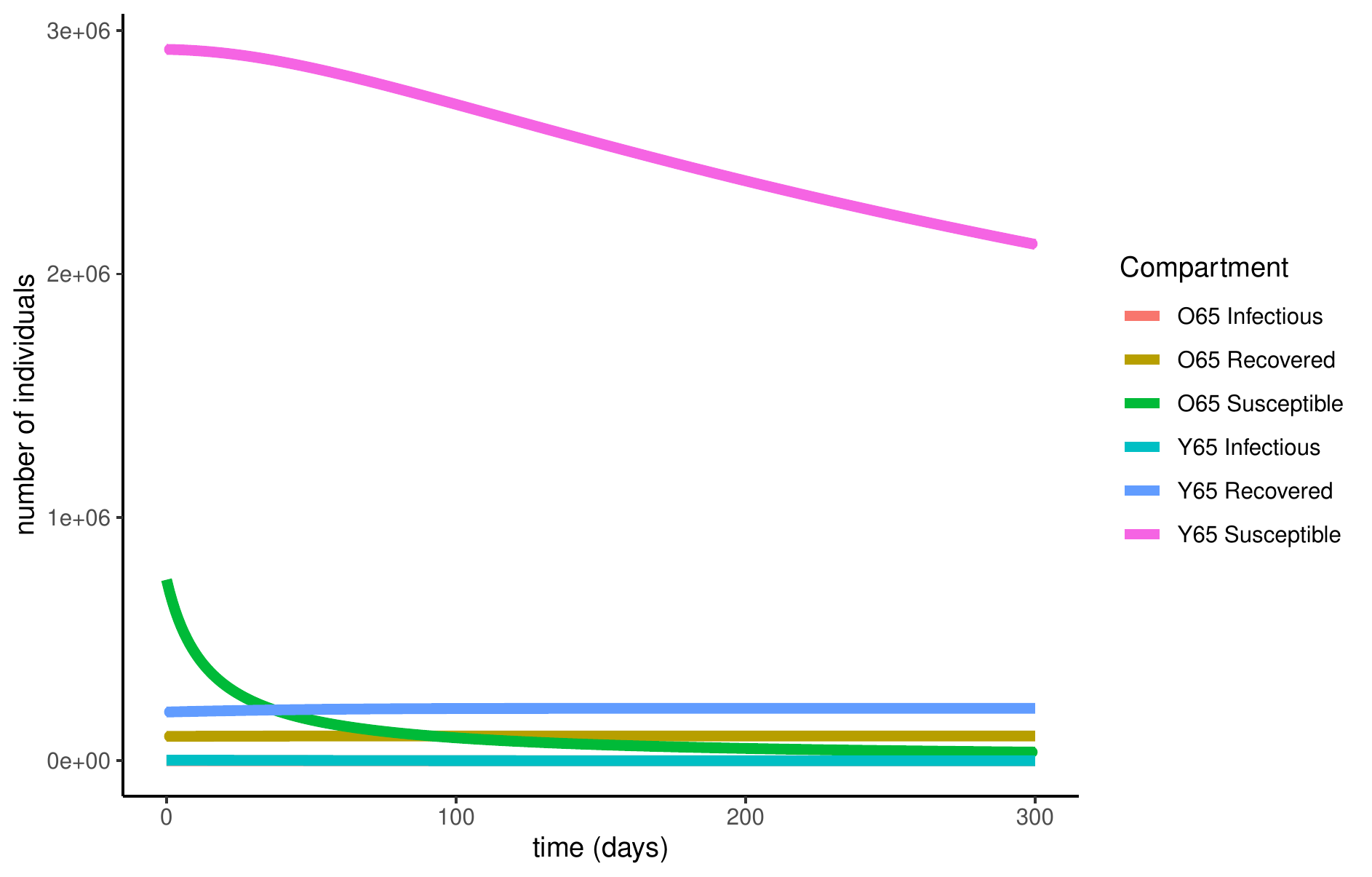}
		\caption{Evolution of states with the vaccine}
	\end{subfigure}
	\caption{ \label{strict measures}Evolution of the pandemic with and without vaccination under strict measures. In the case where no vaccination policy is in place, there is an early peak for the infectious individuals in both age groups (red and turquoise curves). In the second graph the peaks are almost indistinguishable. This is due to the vaccination reducing both the total and the maximum number of infections that take place, hence flattening the curve. The number of susceptible people (pink and blue curves) in both age groups declines fast in the first graph, due to most of them getting infected, while in the second case the decline is much slower, as many individuals are getting vaccinated and thus protected from infection.}
\end{figure}

\begin{figure}[ht]
	\begin{subfigure}{.45\textwidth}
		\includegraphics[width=\linewidth,height=4cm]{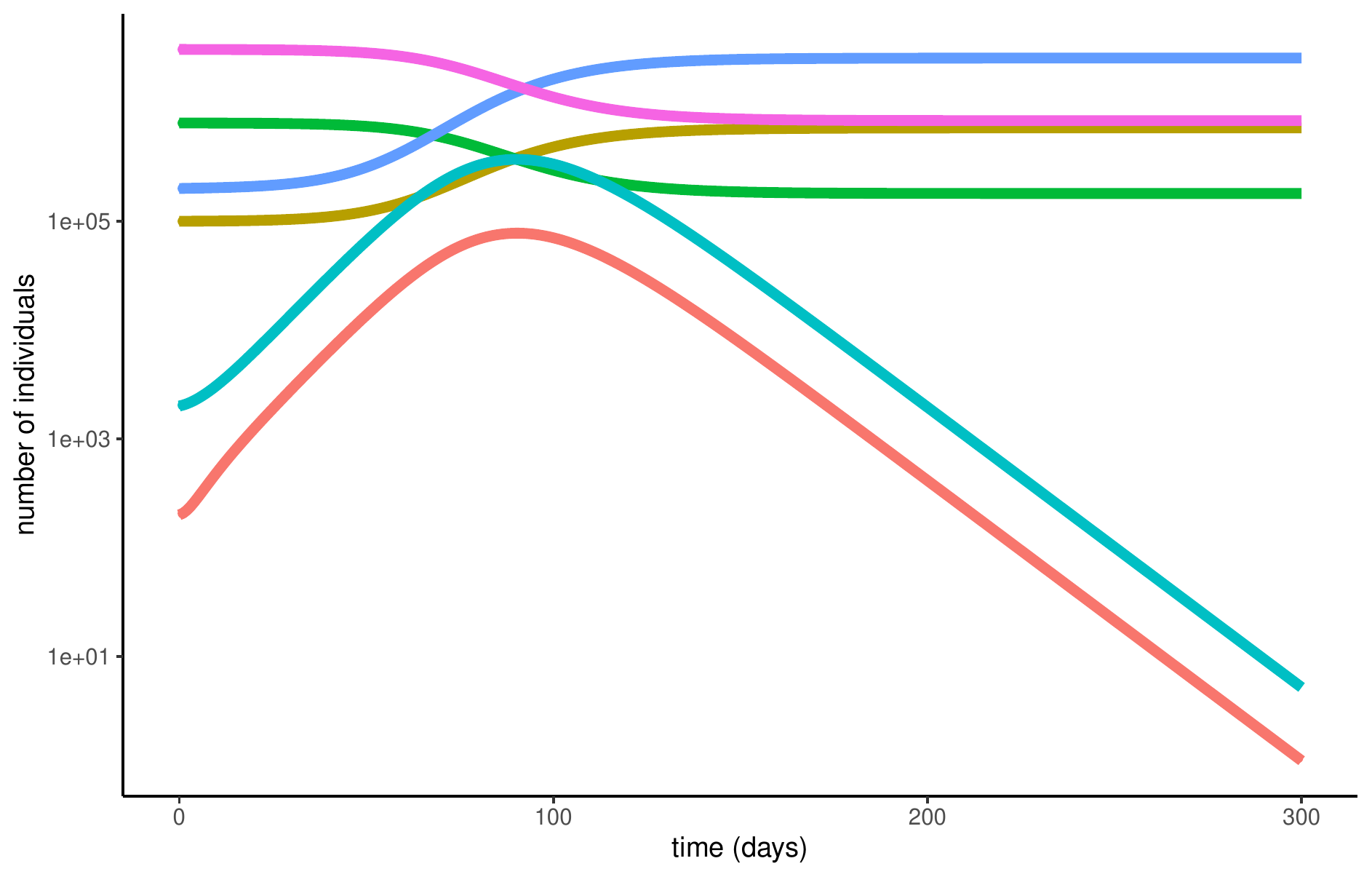}
		\caption{Evolution of states without the vaccine}
	\end{subfigure}%
	\begin{subfigure}{.55\textwidth}
		\includegraphics[width=\linewidth,height=4cm]{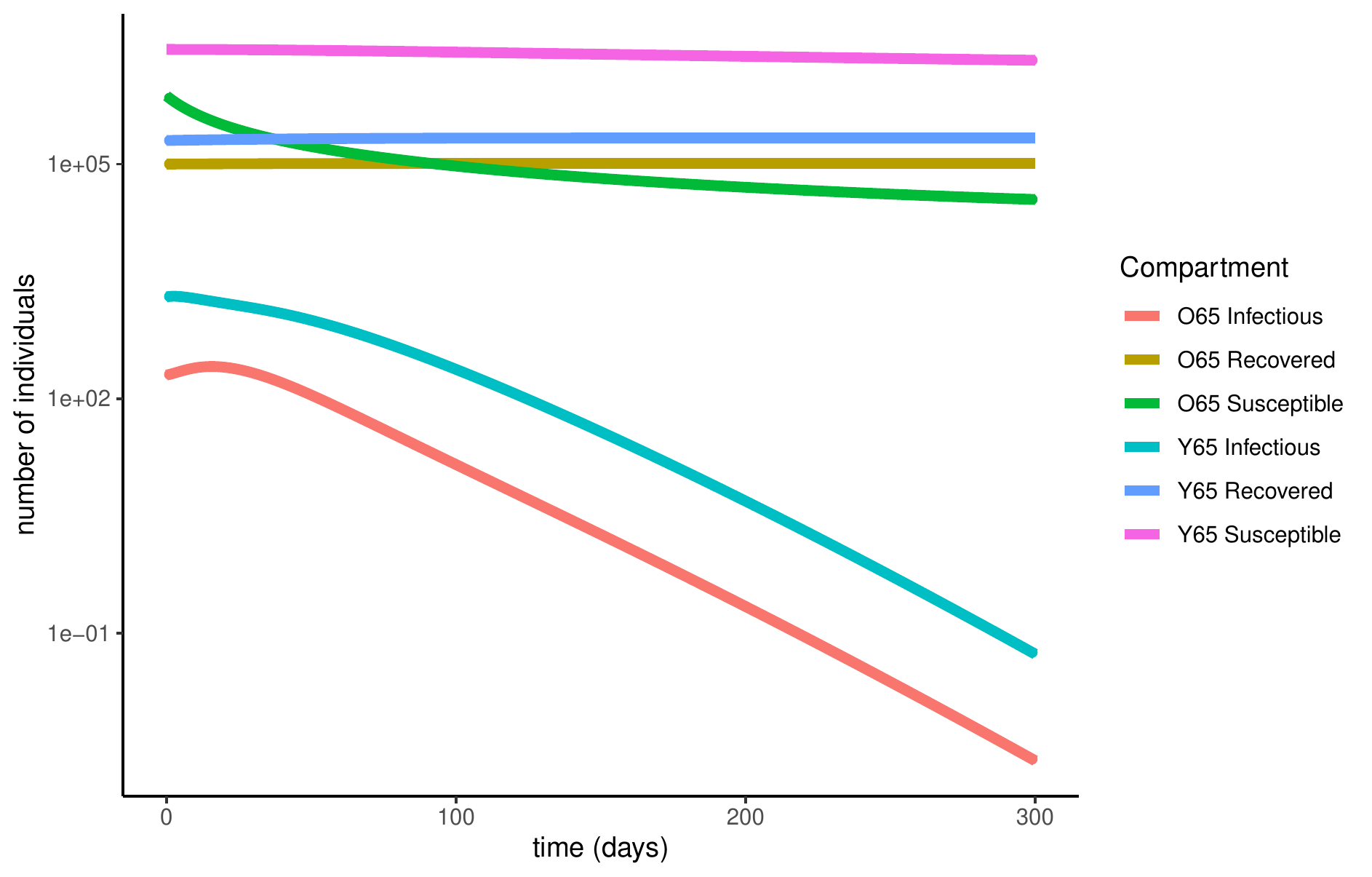}
		\caption{Evolution of states with the vaccine}
	\end{subfigure}
	\caption{\label{strict measures log}Evolution of the pandemic with and without vaccination under strict measures. The graphs are produced in vertical log scale to give a clearer view of the peaks. In the first graph there are clear peaks in the number of infectious individuals in both age groups while in the second, the curves are clearly flattened thanks to the vaccination.}
\end{figure}

\figurename{ \ref{strict measures}} shows the difference between the evolution of the pandemic when there is no control in place and when the optimal control is applied. Due to the transmission rates being quite low, the increase in the number of infectious individuals will not be particularly evident, especially when the vaccination is applied. For that reason we represented the same information also on the log scale in \figurename{ \ref{strict measures log}} to give a clearer picture of the difference an optimal vaccination strategy makes to the evolution of the pandemic. In the case where no control is applied there is a peak in the number of infectious individuals in both age groups and as a result, the number of susceptible people drops fast with the pandemic ending in a very short period of time, namely less than 250 days, with the greatest part of the population having been infected. Specifically, 80\% of the population over 65 (720,249) and 79\% of the population under 65 (3,159,510) gets infected by the virus until the end of the pandemic. That means that around 79\% of the total population of the country (3,879,759) will get infected by the virus when no vaccination strategy is in place. 

On the other hand, when the optimal vaccination strategy is applied, the infectious curve is flattened to the point where it becomes nearly indistinguishable from the horizontal axis. This means that the number of people to become infected by the virus is substantially reduced as a result of the protection provided by the vaccine. In particular, only 11.3\% (101,764) and 5.38\% (215,132) of the over and under 65 years old groups respectively will get infected by the virus, resulting in a total 6.47\% of the population (316,896). This is achieved thanks to the 70.58\% (635,187) and 17.08\% (683,044) of the over 65 and under 65 populations respectively receiving the vaccine and being successfully protected by it. In total, that is 26.9\% of the population (1,318,231) being successfully vaccinated.

\begin{figure}[ht]
	\begin{subfigure}{.45\textwidth}
		\includegraphics[width=\linewidth,height=4cm]{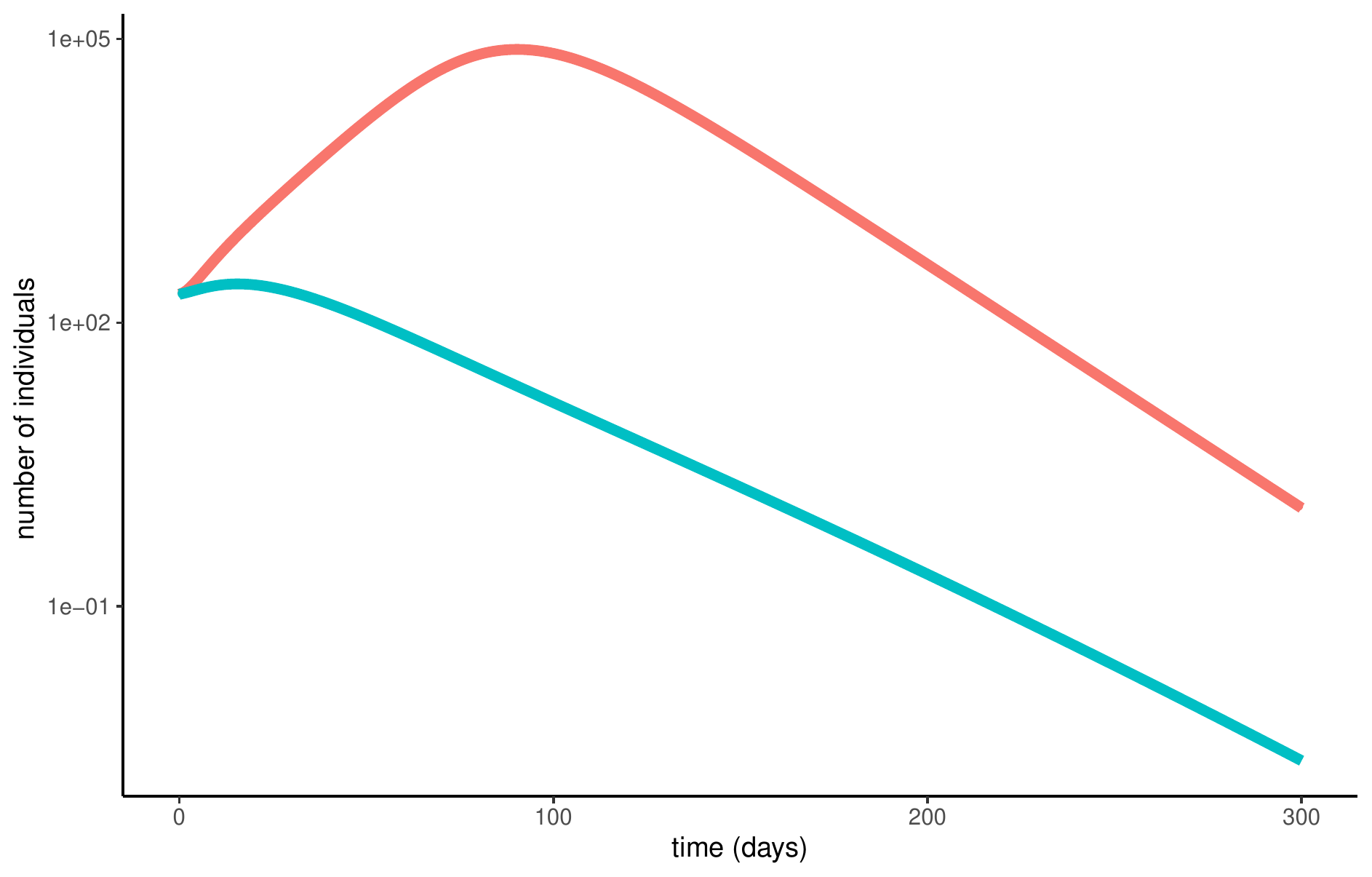}
		\caption{Evolution of infectious in over 65}
	\end{subfigure}%
	\begin{subfigure}{.55\textwidth}
		\includegraphics[width=\linewidth,height=4cm]{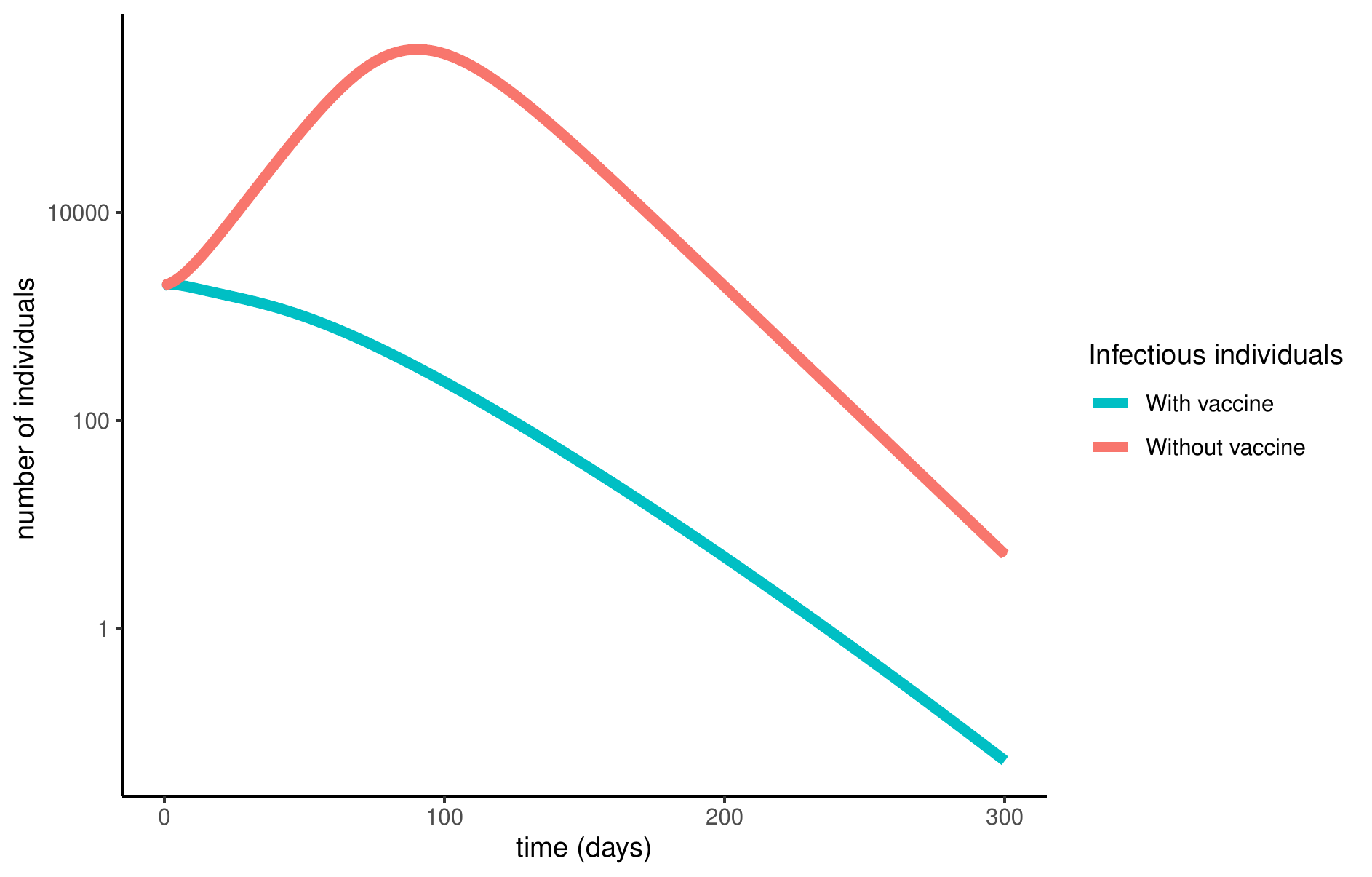}
		\caption{Evolution of infectious in under 65}
	\end{subfigure}
	\caption{\label{strict measures infectious}Comparison of baseline curves (without vaccine) with optimal vaccination strategy curves for the infectious populations in both age groups under strict measures. The vertical scale is in  log units. The vaccination successfully flattens the curve in both cases, resulting in fewer infections overall as well as fewer simultaneous infections. }
\end{figure}

The difference that the control function makes to the number of infectious individuals (in the log scale) can be more clearly seen in \figurename{ \ref{strict measures infectious}}. For both age groups, the curve is substantially flattened as a result of the vaccination. This is a visual representation of the reduction in the total number of infections and also showcases the fact that less people will be infected at the same time, causing less stress on the healthcare system.    

\begin{figure}[ht]
	\begin{subfigure}{.32\textwidth}
		\includegraphics[width=\linewidth,height=3.5cm]{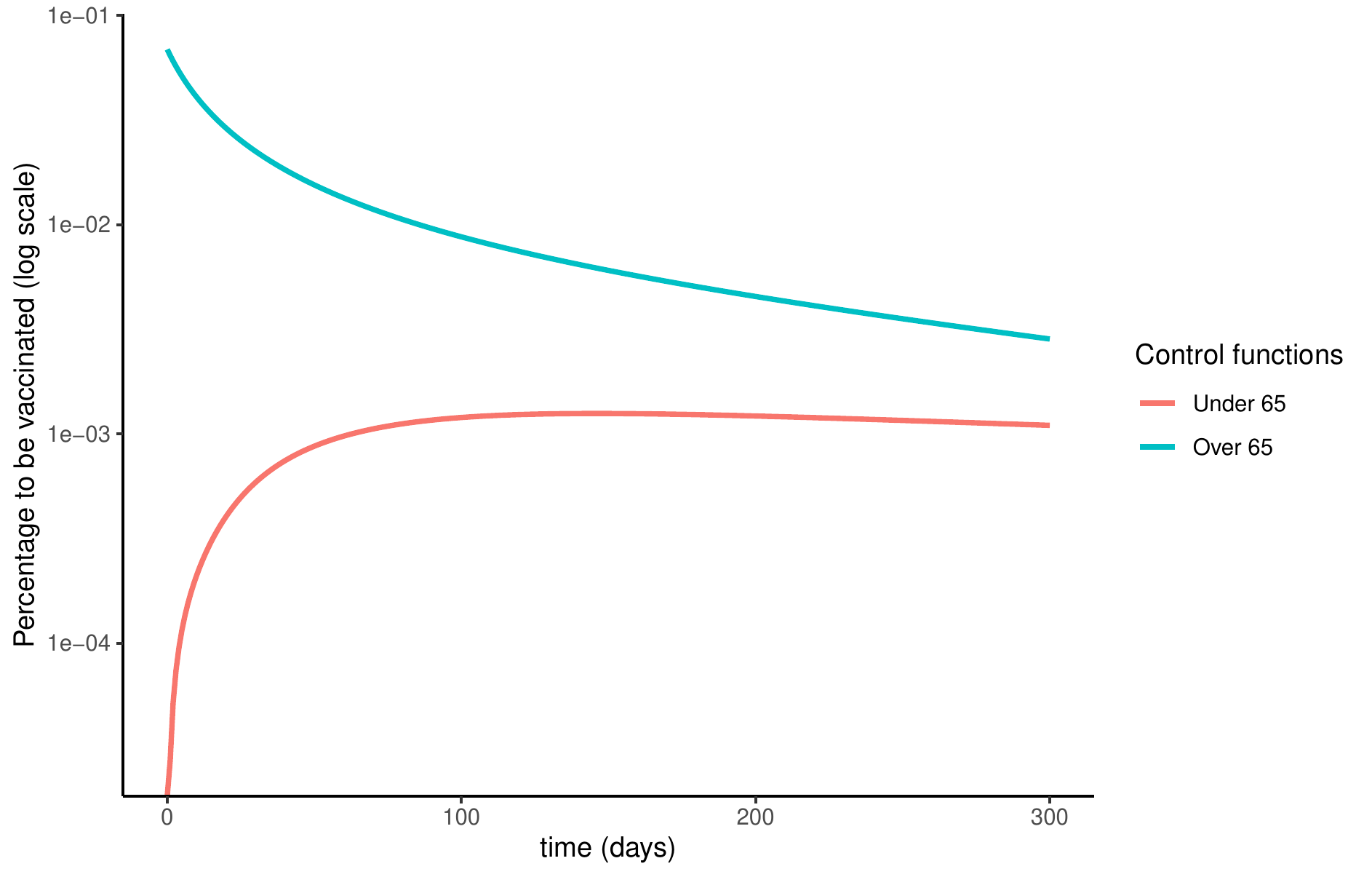}
		\caption{\label{strict measures controls log}Optimal control functions}
	\end{subfigure}%
	\begin{subfigure}{.32\textwidth}
		\includegraphics[width=\linewidth,height=3.5cm]{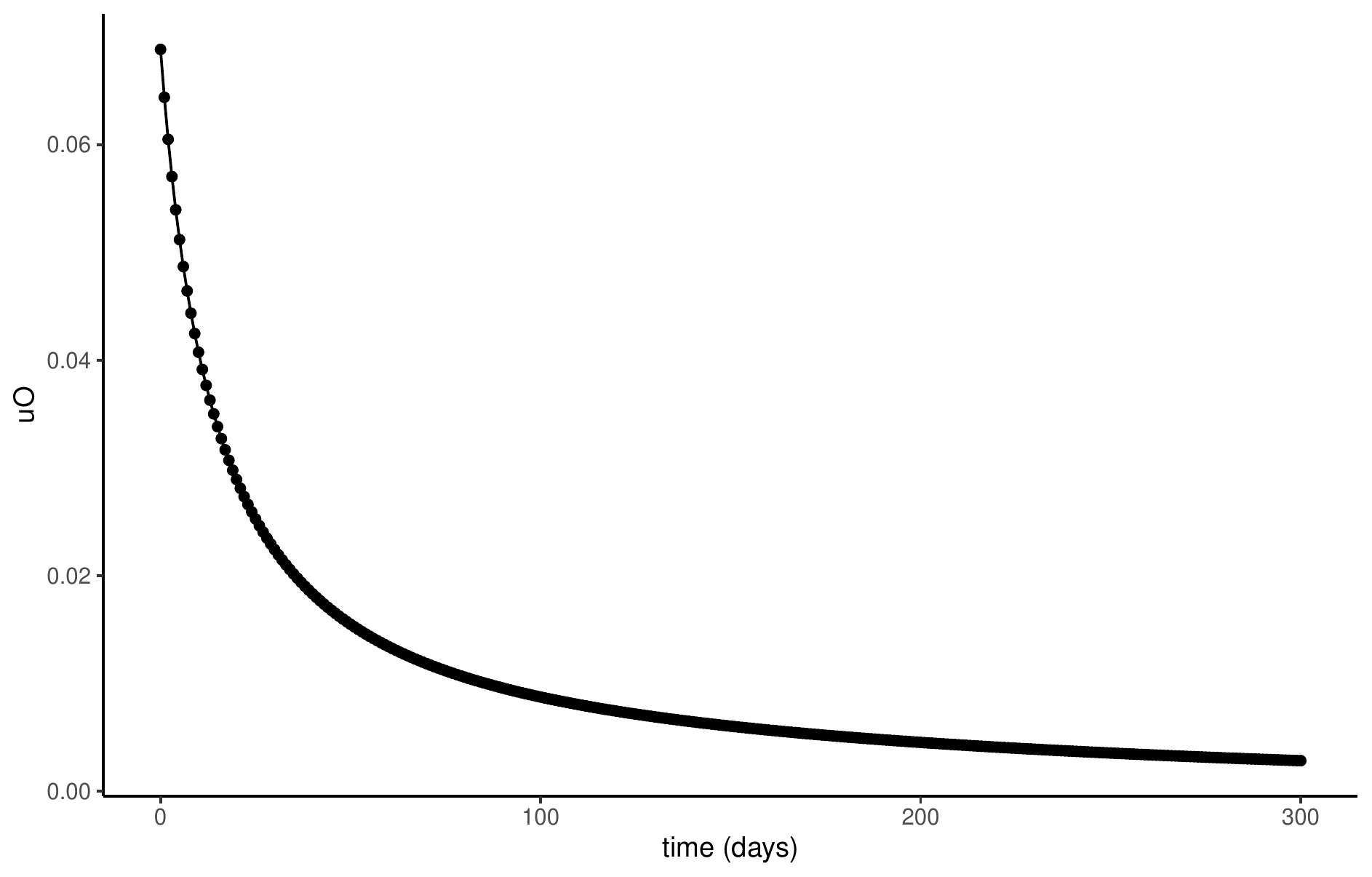}
		\caption{Optimal control for over65}
	\end{subfigure}
	\begin{subfigure}{.32\textwidth}
		\includegraphics[width=\linewidth,height=3.5cm]{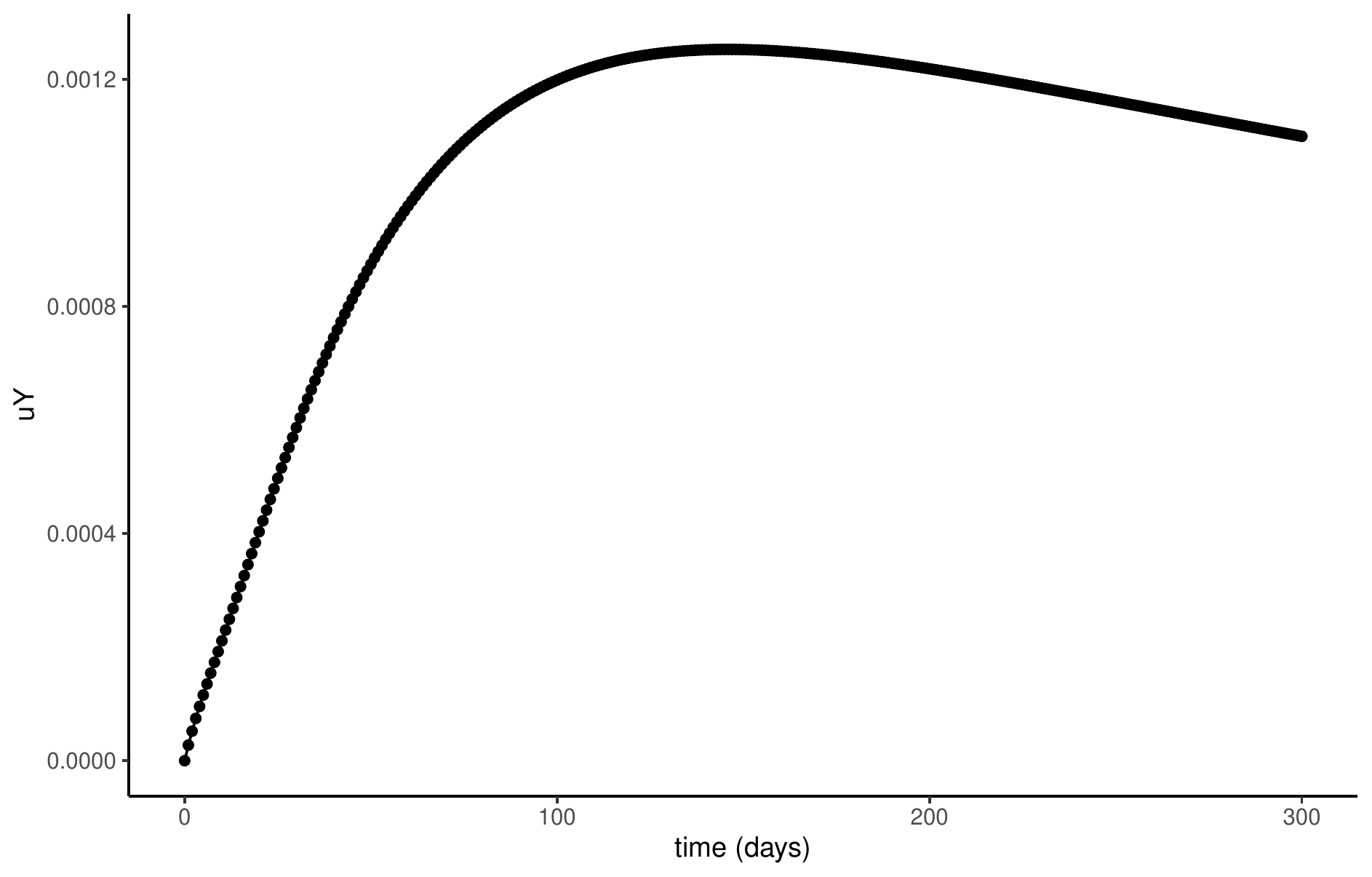}
		\caption{Optimal control for under65}
	\end{subfigure}
	\caption{\label{strict measures controls} Optimal control functions for both age groups. \ref{strict measures controls log} presents both control functions on a vertical log scale. The o65 group gets vaccinated faster and in higher percentages, while the vaccination of the y65 group commences but proceeds slower and continues after the pandemic is under control.}
\end{figure}

The optimal vaccination strategy for each age group, i.e., the control functions that resulted from the optimisation technique described in \sectionautorefname{ \ref{section optimisation}} can be seen in \figurename{ \ref{strict measures controls}}. As expected, close to the start of our simulations, a high percentage of the susceptible population gets vaccinated to quickly get the pandemic under control, and then that percentage continuously declines. It is interesting to note how the o65 population gets vaccinated earlier and with higher percentages while the vaccination of the younger population is happening slower, with smaller percentages every day, and goes on after the pandemic is over as can be seen in \figurename{ \ref{strict measures controls}}. This is also evident in the percentages of the people in each age group that get vaccinated by the end of the simulation, with 70.58\% and 17.08\% for the over and under 65 groups respectively, as described earlier. 

It is worth noting that the parameter setup used for this simulation is indicative of a national lockdown situation, which is an unrealistic scenario for this length of time. For this reason, we looked into a second parameter setup, discussed in the next section.

\subsection{\label{section loose restrictions}Evolution of pandemic under loose restrictions}

We now simulate the evolution of the pandemic in Ireland, using the same initial conditions for the states as defined in \tableautorefname{ \ref{initial states}}. We are interested in investigating the effect a different choice of parameters, specifically different $R0$ numbers, would have on our results. We are looking into values that are indicative of an extreme scenario where the measures in place are not strict or even non existent and each infected individual can infect multiple individuals before recovering. The column titled Case 2 in \tableautorefname{ \ref{initial parameters}} consists of values chosen for the parameters, while the rest of the parameters remain the same as in \subsectionautorefname{ \ref{section tight restrictions}}. 

\begin{figure}[!ht]
	\begin{subfigure}{.45\textwidth}
		\includegraphics[width=\linewidth,height=4cm]{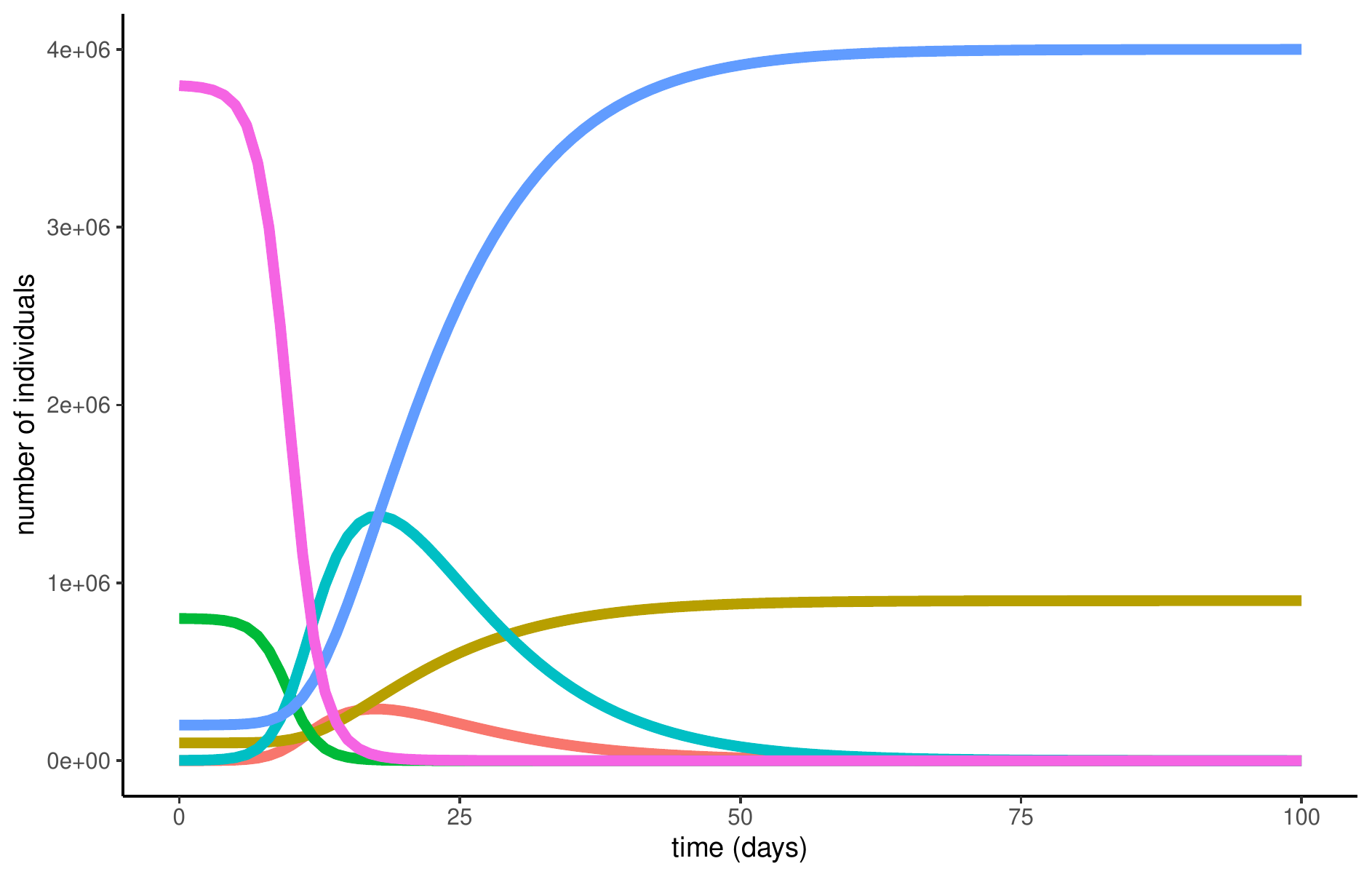}
		\caption{Evolution of states without the vaccine}
	\end{subfigure}%
	\begin{subfigure}{.55\textwidth}
		\includegraphics[width=\linewidth,height=4cm]{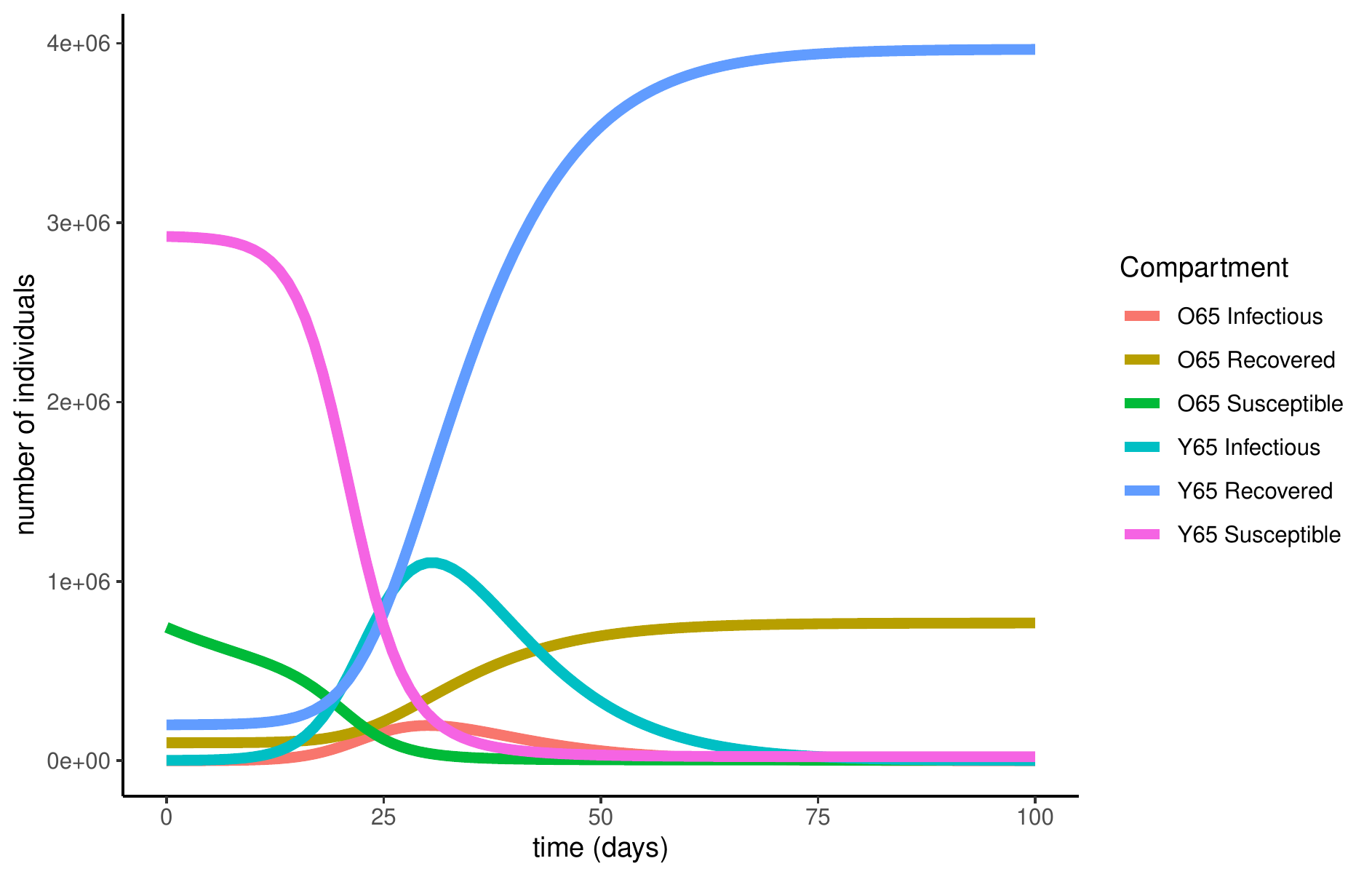}
		\caption{Evolution of states with the vaccine}
	\end{subfigure}
	\caption{\label{loose restrictions}Evolution of the pandemic with an without vaccination under loose measures. Without a vaccination policy in place, the curve of infectious individuals presents an earlier and higher peak as opposed to the case where the vaccination strategy is applied and the curve is flattened. }
\end{figure}

\begin{figure}[ht]
	\begin{subfigure}{.45\textwidth}
		\includegraphics[width=\linewidth,height=4cm]{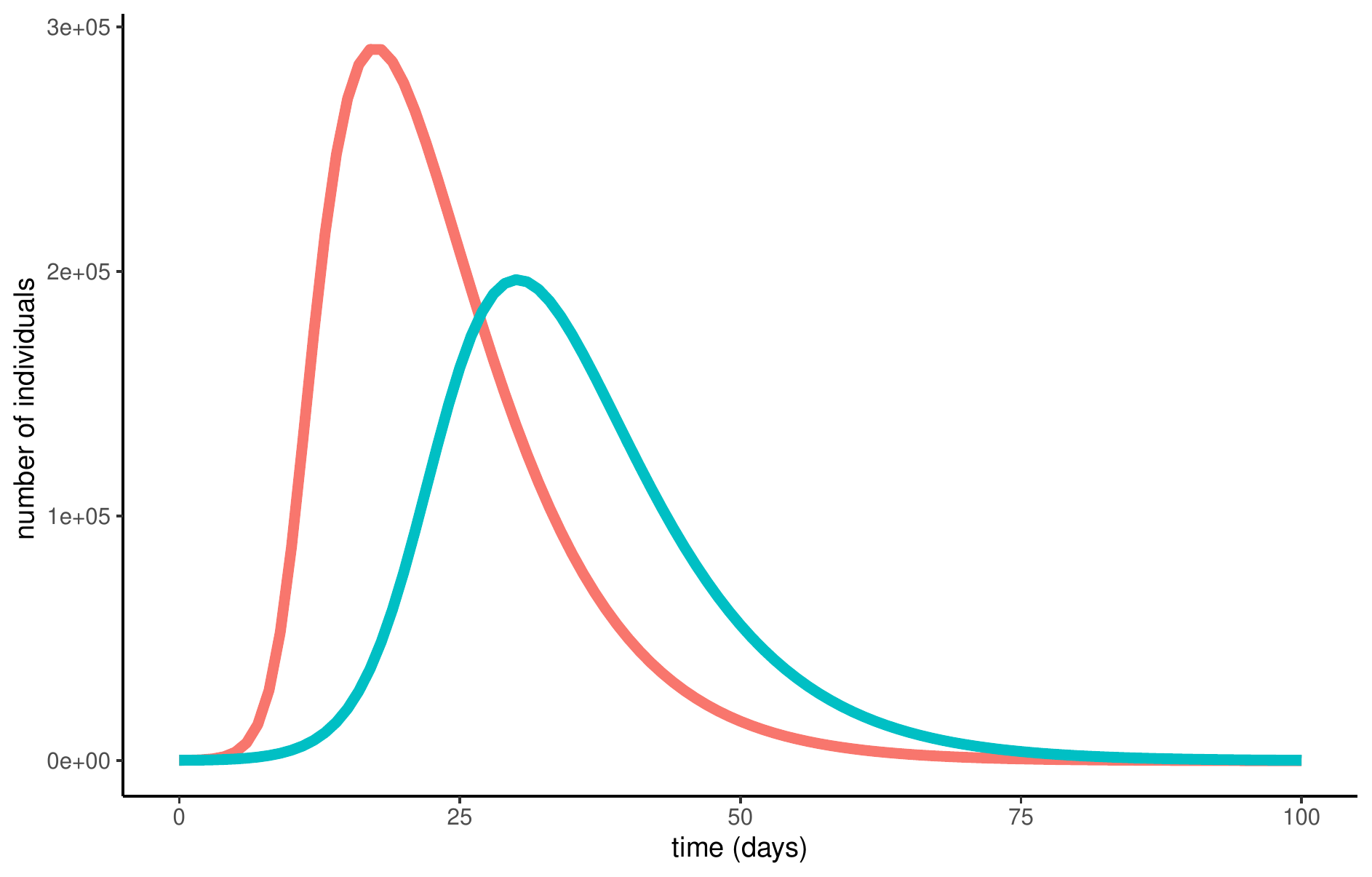}
		\caption{Evolution of infectious in over 65}
	\end{subfigure}%
	\begin{subfigure}{.55\textwidth}
		\includegraphics[width=\linewidth,height=4cm]{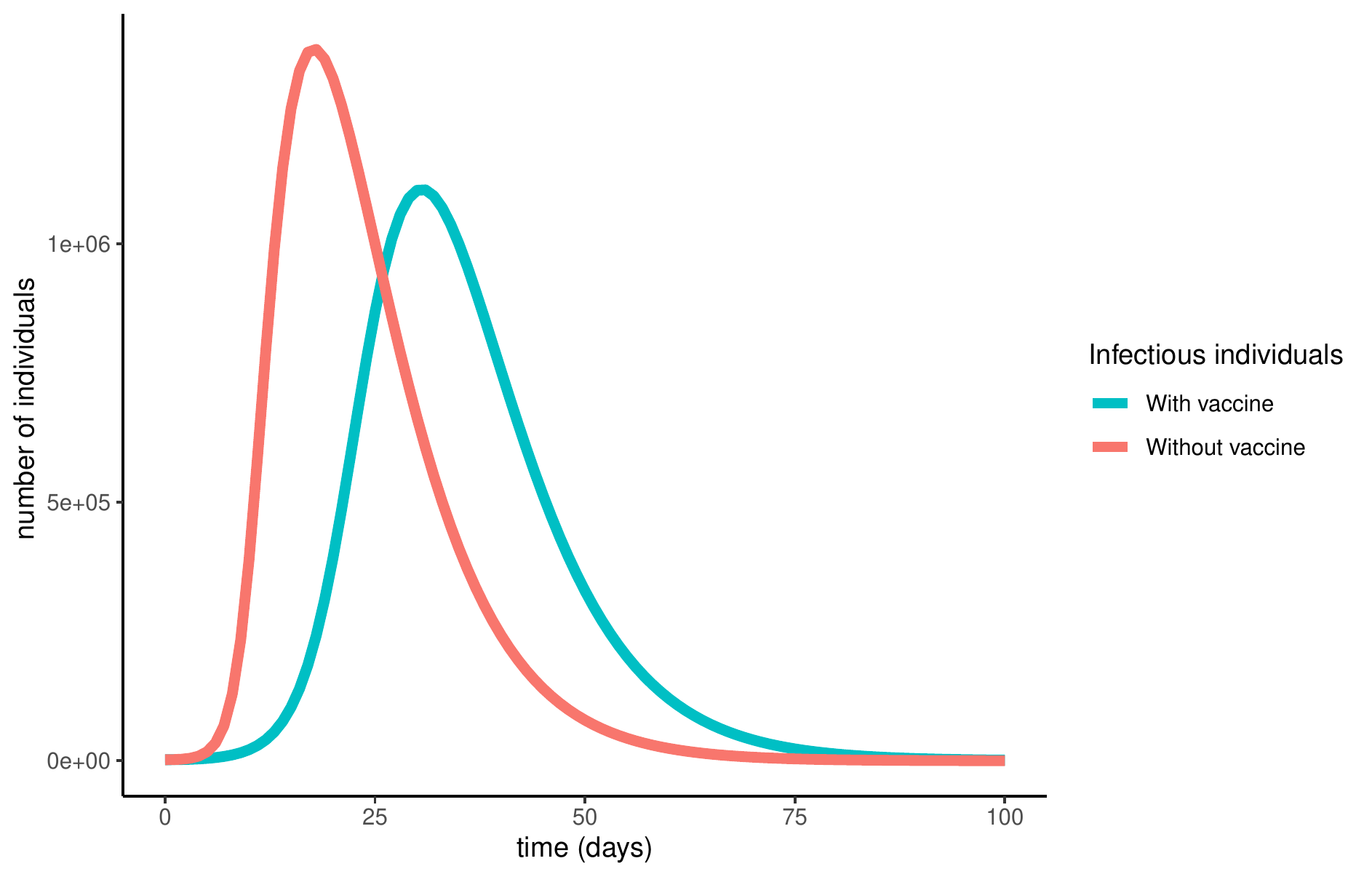}
		\caption{Evolution of infectious in under 65}
	\end{subfigure}
	\caption{\label{loose restrictions infectious}Comparison of baseline curves (without vaccine) with optimal vaccination strategy curves for the infectious populations in both age groups under loose restrictions. }
\end{figure}

When no vaccination policy is in place, the curves of infectious individuals in both age groups present a very high peak due to the high transmission rates as can be seen in Figures \ref{loose restrictions} and \ref{loose restrictions infectious}. This peak is indicative of a period of time where a very large number of people are infected at the same time, a situation that would cause a huge strain on the country's healthcare system. The percentages of people that get infected in the over and under 65 group respectively are 99.995\% (899,954) and 99.993\% (3,999,713), resulting in a total 99.993\% (4,899,667) of the population. This means that nearly every single resident of the country will get infected by the virus in a very short period of time, while at the same time paralysing the healthcare system. 

That curve gets flattened, with a lower peak and greater spread, thanks to the optimal control functions, i.e. the optimal vaccination strategy. This would offer some relief to the healthcare system, as not all infections will be active at the same time. The percentages of people that get infected in the over and under 65 group respectively are 85.29\% (767,614) and 99.1\% (3,965,963), resulting in a total 96.6\% of the population (4,733,577). These percentages, though reduced, are still very high and that is due to the fact that the chosen transmission rates are so high. For the same reason the percentages of the people that get successfully protected from the vaccine are relatively low, namely 14.2\% (127,871) and 0.1\% (3,861) for the over 65 and under 65 groups respectively, meaning a total 2.7\% of the population (131,732) will be successfully vaccinated. There is an additional reason to why these percentages are so low. From the moment an individual receives the vaccine, there is a period equal to the \textit{mean holding time vaccinated} ($t_V$) during which the person may still get infected. Due to the high transmission rates, a lot of vaccinated individuals will get infected, hence not making it to the protected state.

A side by side presentation of the evolution of the pandemic with and without vaccination is given in \figurename{ \ref{loose restrictions}} and a comparison of the optimised with the baseline curves for the infectious populations in each age group are given in \figurename{ \ref{loose restrictions infectious}}.

\begin{figure}[ht]
	\begin{subfigure}{.32\textwidth}
		\includegraphics[width=\linewidth,height=3.5cm]{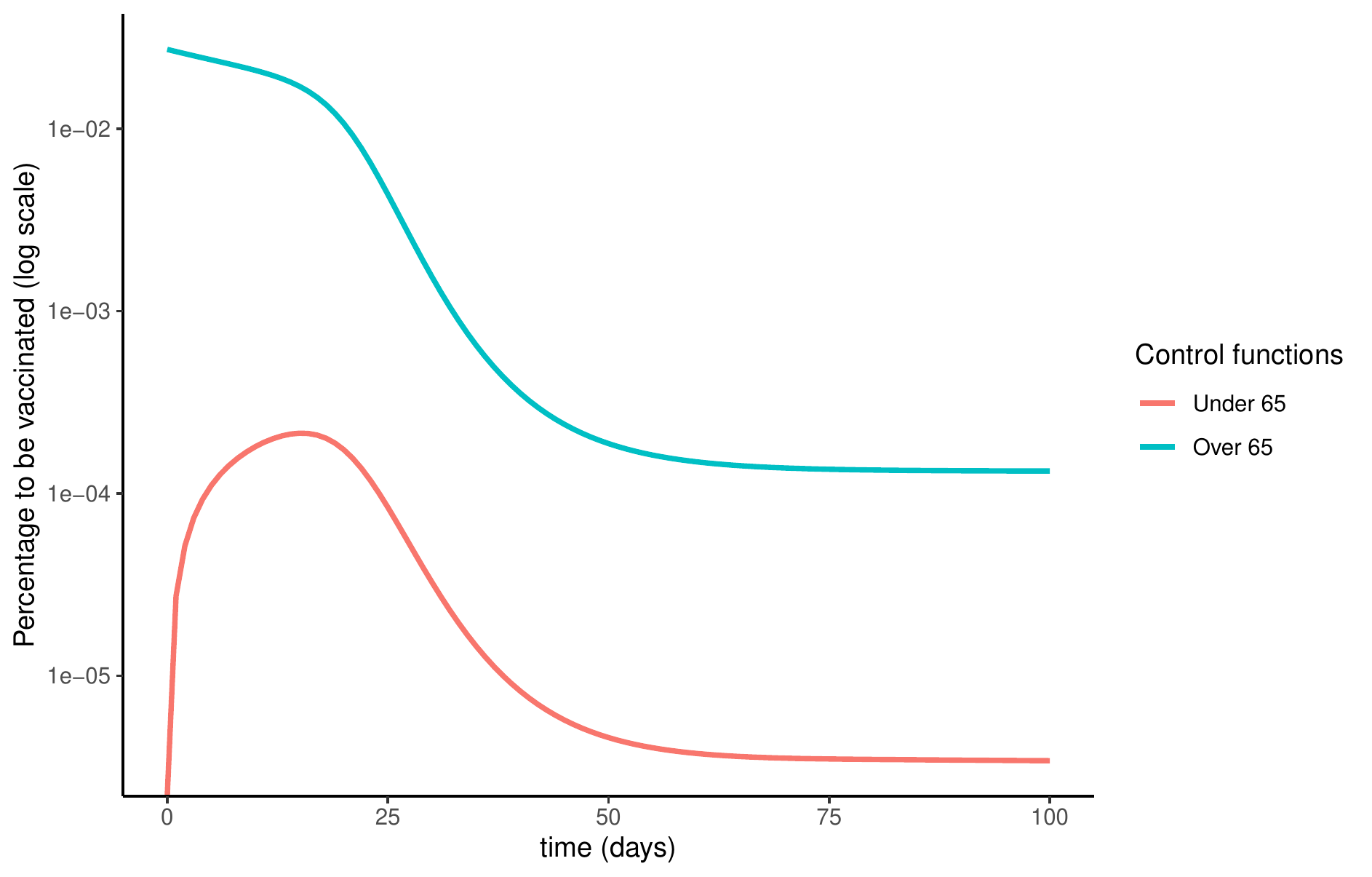}
		\caption{\label{loose restrictions optimal control log}Optimal control functions}
	\end{subfigure}%
	\begin{subfigure}{.32\textwidth}
		\includegraphics[width=\linewidth,height=3.5cm]{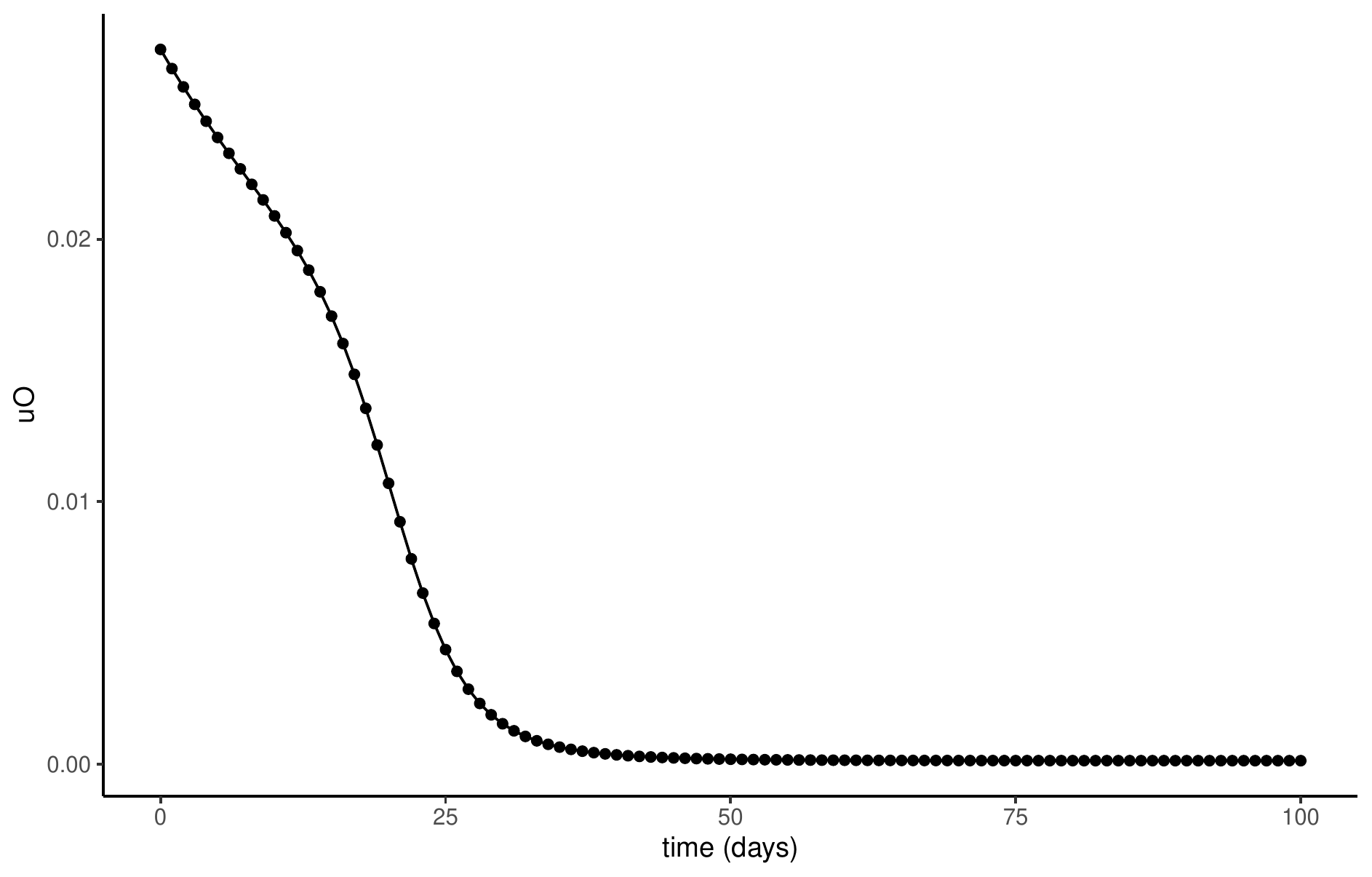}
		\caption{Optimal control for over65}
	\end{subfigure}
	\begin{subfigure}{.32\textwidth}
		\includegraphics[width=\linewidth,height=3.5cm]{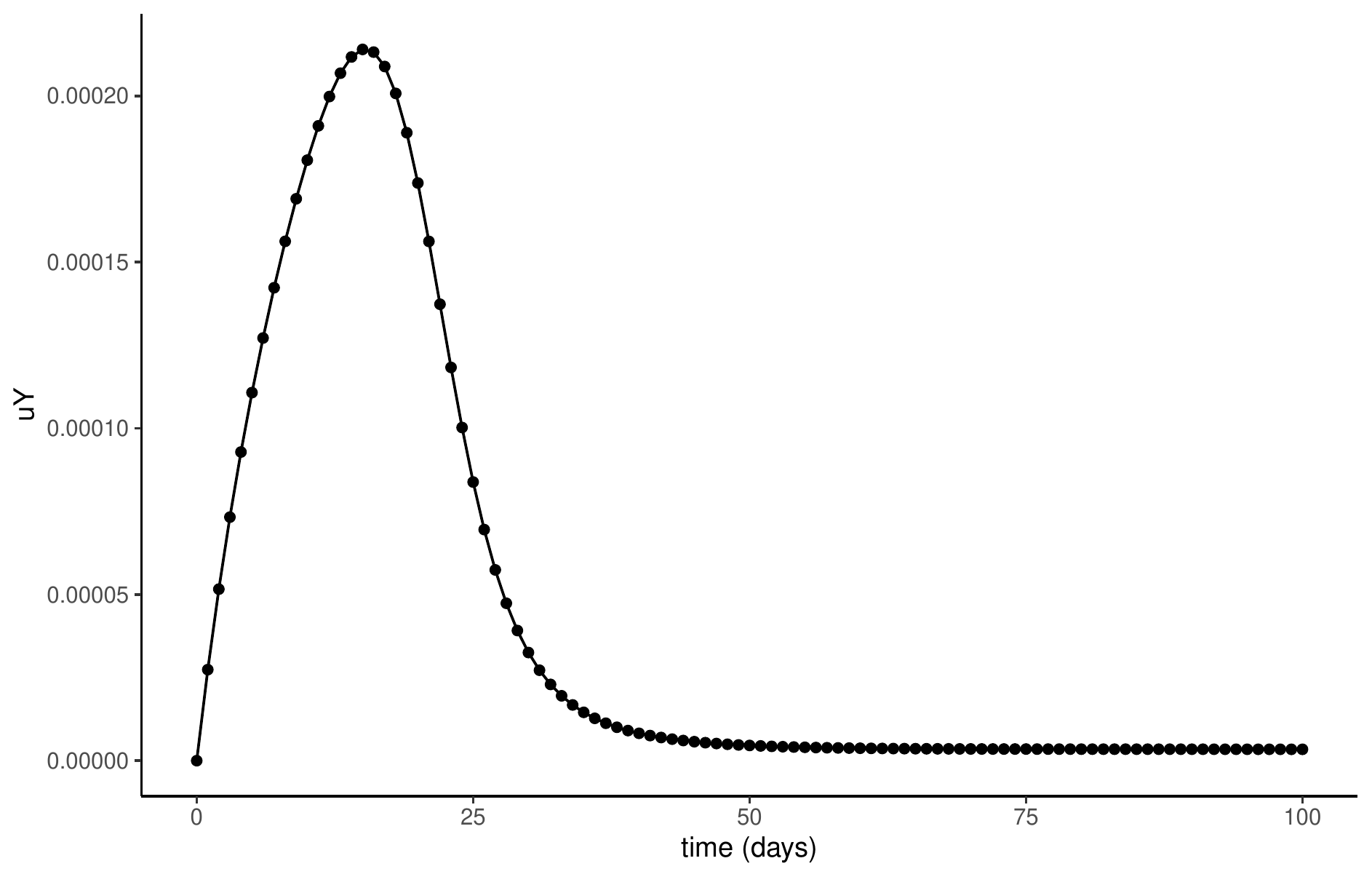}
		\caption{Optimal control for under65}
	\end{subfigure}
	\caption{\label{loose restrictions optimal controls}Optimal control functions for both age groups. \ref{loose restrictions optimal control log} presents both control functions in log scale.}
\end{figure}

\figurename{ \ref{loose restrictions optimal controls}} shows the optimal control functions that can be obtained as the solution to our optimisation problem. Similar to the previous simulations, the control function curves present a peak at first, indicating a high percentage of the population being chosen for vaccination early on and that percentage getting reduced as time progresses.

\section{Comparative Study}\label{section comparison}

\begin{figure}[!ht] 
	\begin{subfigure}{.5\textwidth}
		\includegraphics[width=\linewidth,height=4cm]{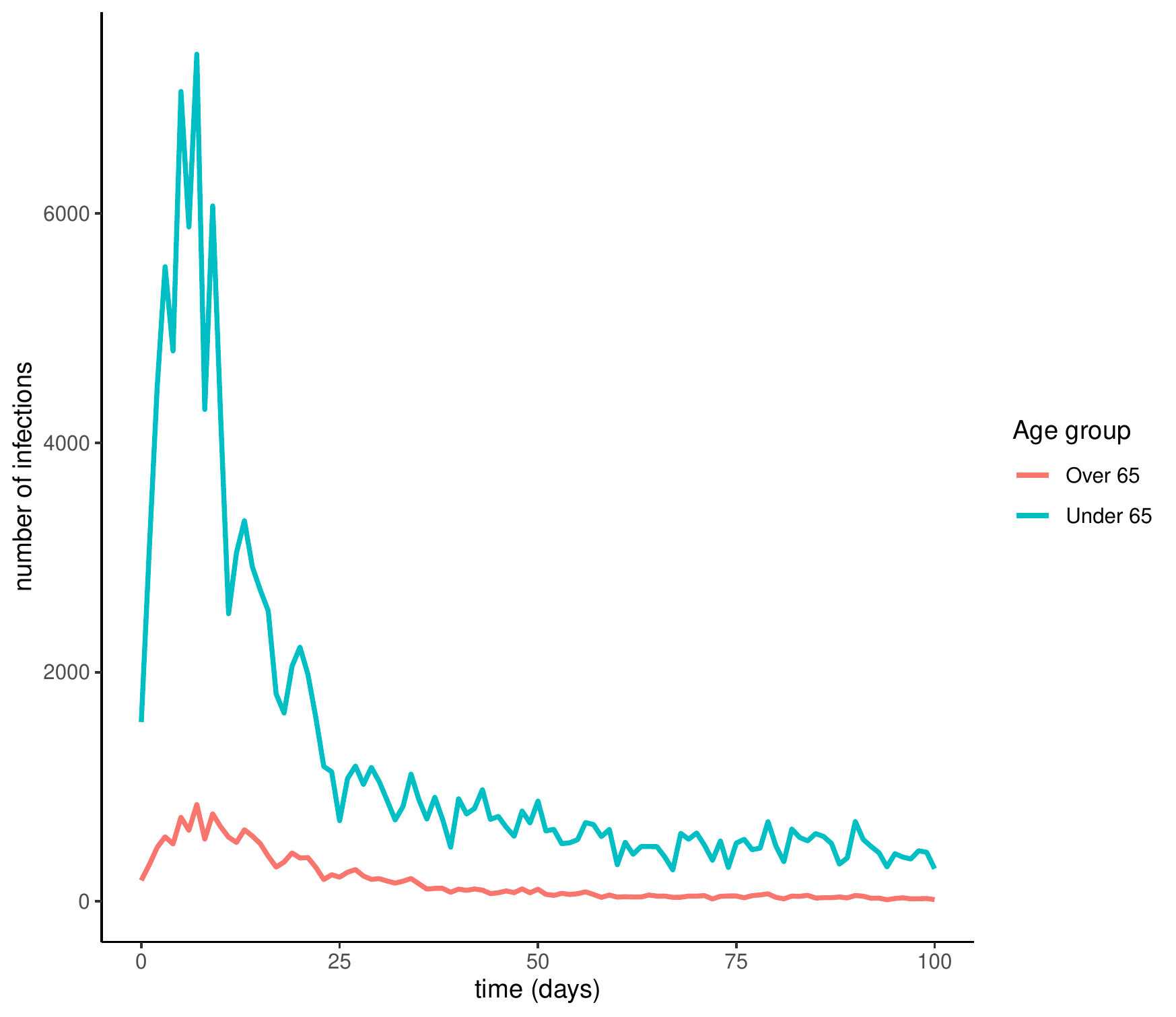}
		\caption{Infectious individuals}
	\end{subfigure}%
	\begin{subfigure}{.5\textwidth}
		\includegraphics[width=\linewidth,height=4cm]{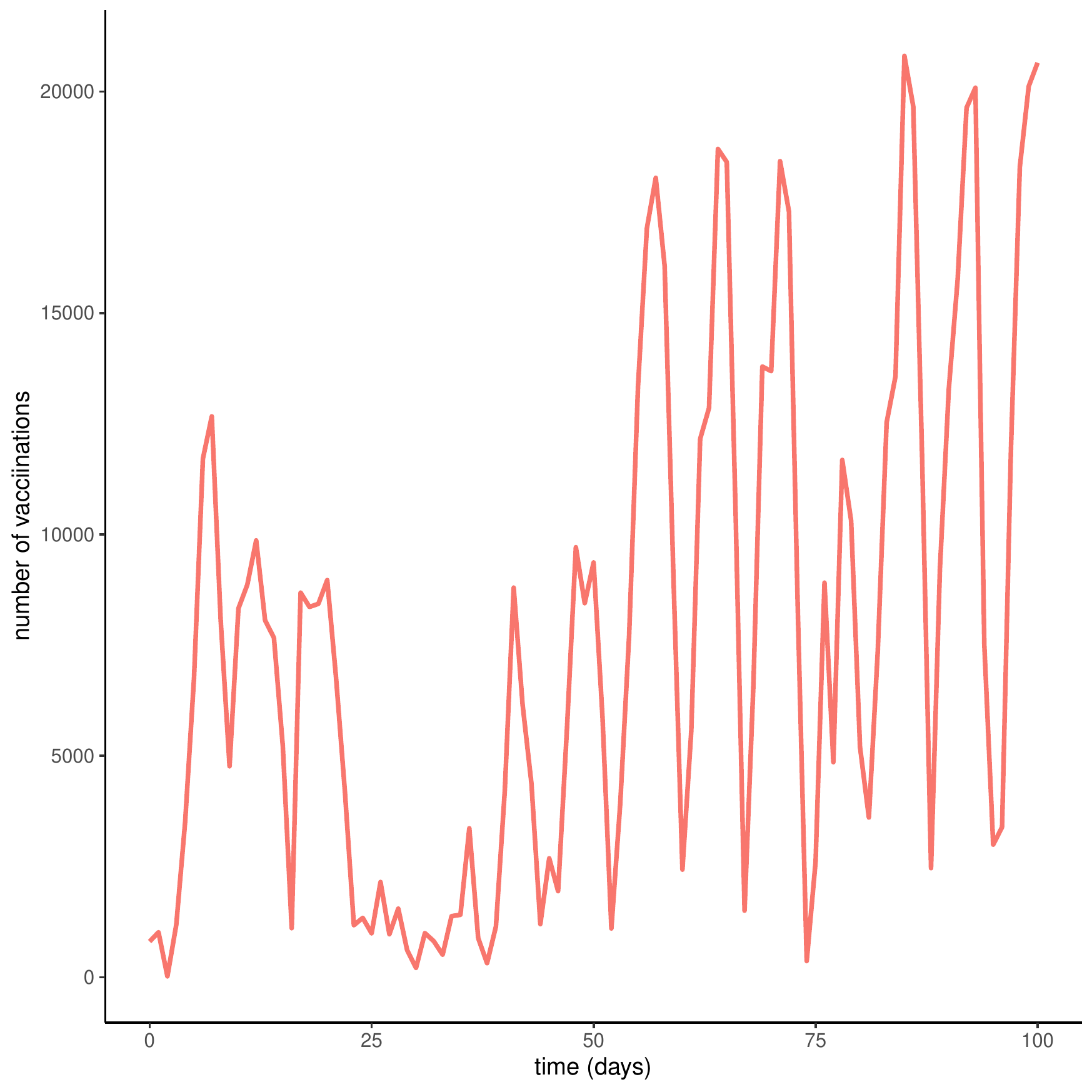}
		\caption{\label{vac before peak}Vaccines administered}
	\end{subfigure}
	\caption{\label{inf and vac before peak}Infection and vaccination numbers for the Republic of Ireland starting from January 
		1\textsuperscript{st} 2021}
\end{figure}

We now look at the effect that the optimal control obtained from our model would have on the evolution of the pandemic in Ireland by comparing our approach to the real infection numbers, given the vaccination strategy taken by the Irish government. In order for the comparison to be sensible, we limit our simulations to a time frame of 100 days, thus ensuring that the transmission rates remain somewhat constant. This is needed because our model does not account for varying $R0$ numbers, while in reality these numbers can fluctuate depending on the restrictive measures in place (e.g. lockdowns) and the preventative measures taken, namely vaccinations. We obtained our data from \cite{datahub}, specifically the daily infection and vaccination numbers starting from January 1\textsuperscript{st} 2021. These values are shown in \figurename{\ref{inf and vac before peak}}.

\begin{figure}[!ht] 
	\begin{subfigure}{.5\textwidth}
		\includegraphics[width=\linewidth,height=4cm]{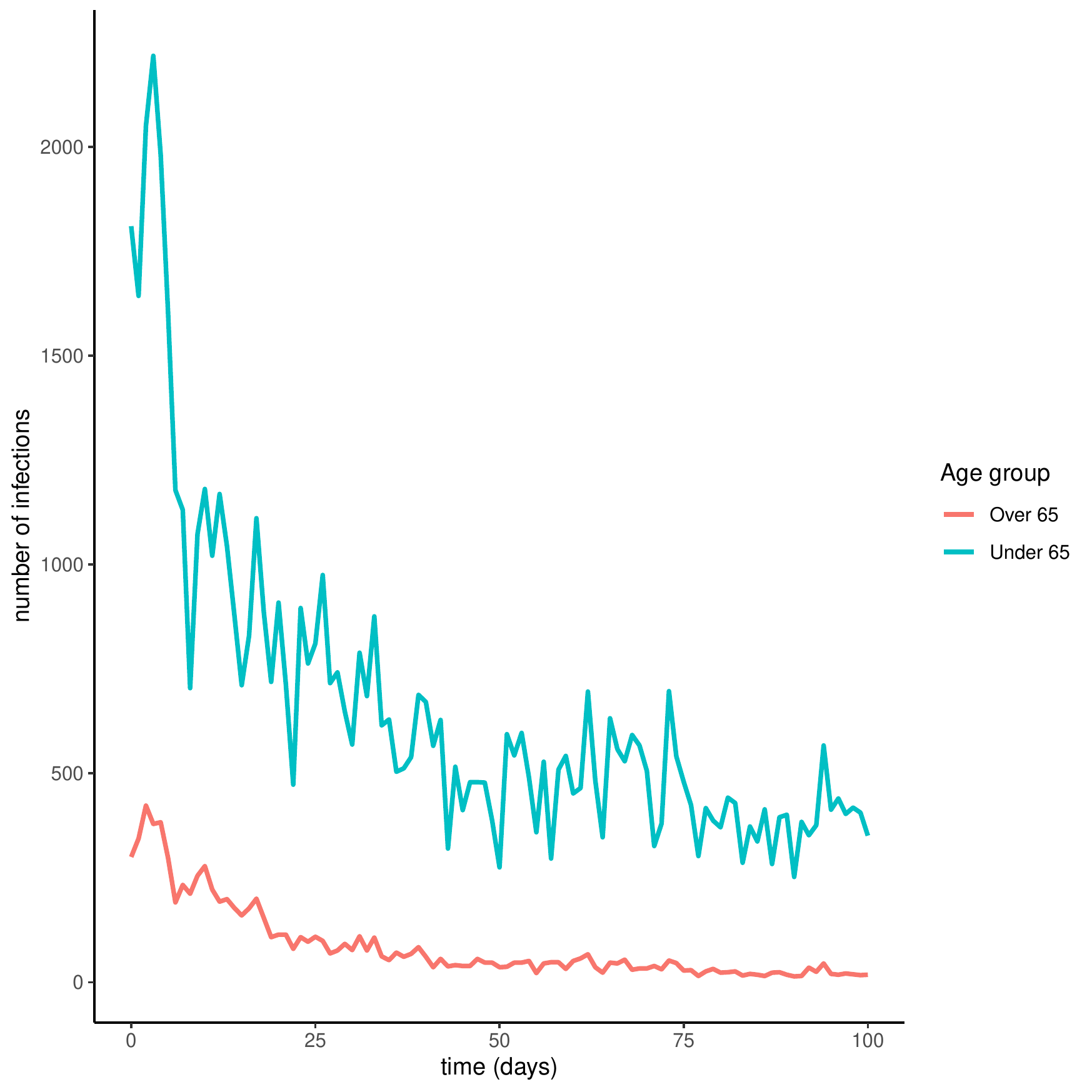}
		\caption{\label{inf after peak}Infectious individuals}
	\end{subfigure}%
	\begin{subfigure}{.5\textwidth}
		\includegraphics[width=\linewidth,height=4cm]{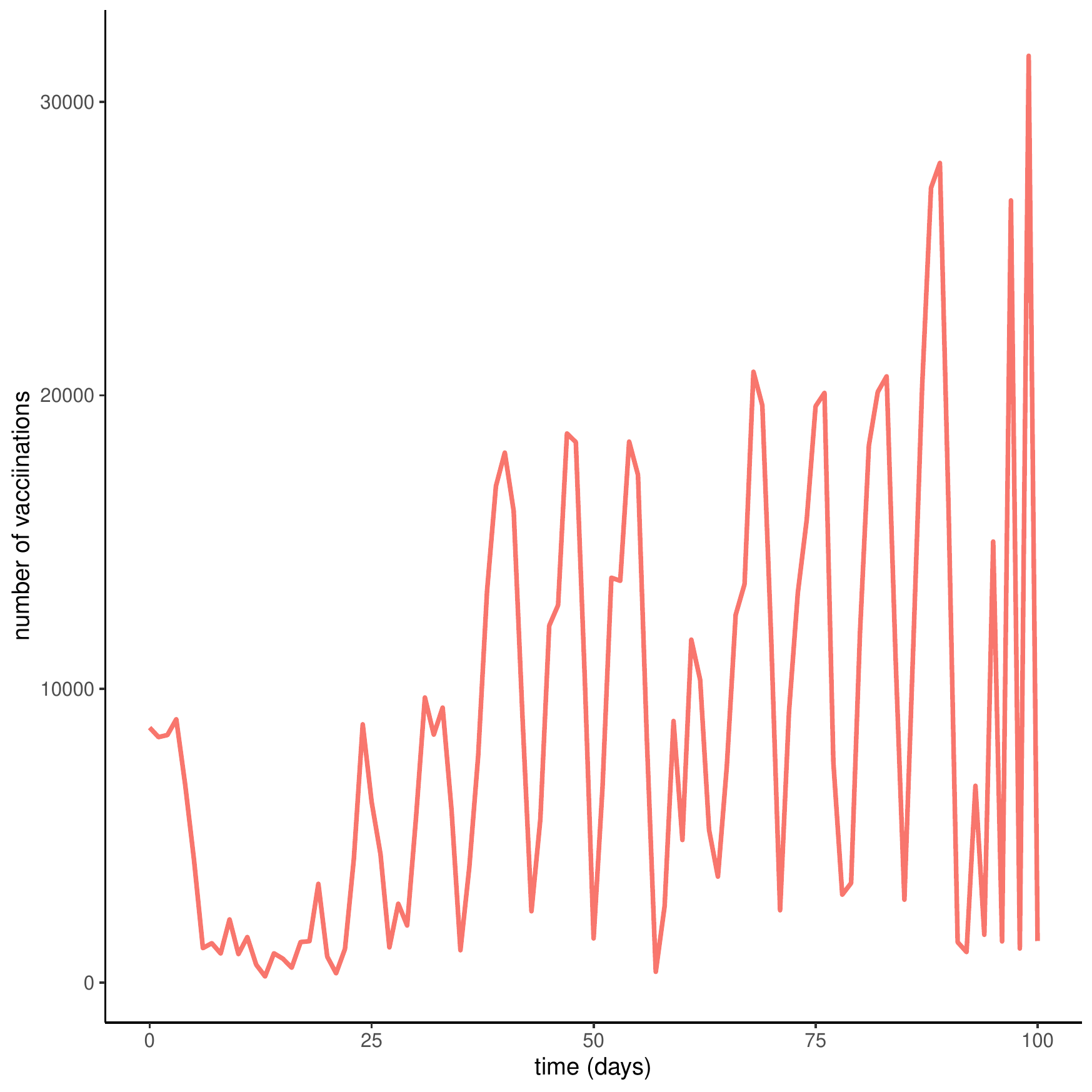}
		\caption{\label{vac after peak}Vaccines administered}
	\end{subfigure}
	\caption{\label{inf and vac after peak}Infection and vaccination numbers for the Republic of Ireland starting from January 
		18\textsuperscript{th} 2021}
\end{figure}

It is evident that there is a substantial peak in the infection numbers during the first few days which can be attributed to the looser measures in place before the Christmas holidays of 2020. Considering this, we produce a second set of comparative plots, starting the simulation just after the peak in infections, namely on January 18\textsuperscript{th} instead of the 1\textsuperscript{st}. The daily infection and vaccination numbers for a period of 100 days starting from January 18\textsuperscript{th} can be seen in \figurename{\ref{inf and vac after peak}}.

Our model requires a specific upper bound for the control functions in order for the optimisation objective to be defined, as seen in Equation \ref{minimization of functional}. One choice for that upper bound could be the average of the daily vaccinations over the time period we are studying. However, it is evident from Figures \ref{vac before peak} and \ref{vac after peak} that there is a great fluctuation in the number of vaccines that were administered which can be attributed to the limited availability of vaccines at the start of 2021. This leads us to assume that the number of vaccines administered each day was the maximum of available vaccines for that day. With that in mind, we designed a time varying upper bound, using each daily vaccination number as the upper bound for the control function for that day. Additionally, while the data provides the total number of vaccines to be administered to the population each day, our model accounts for two age groups, each of them requiring a separate upper bound. To address this issue, we apply complimentary percentages of the daily vaccination numbers as upper bounds to each control function. Specifically, and in agreement with the policy taken by the Irish government, we consider a higher percentage of the vaccines to be made a available to the over 65 group. We choose 80\% of the available vaccinations to be used for the over-65 group and let $V_{No}(t)$ be the number of vaccines administered at time (day) $t$. Then, for time $t$, the control functions are bounded as follows:
\begin{equation*}
u_O(t)S_O(t)<0.8 \cdot V_{No}(t) \text{ and } u_Y(t)S_Y(t)<0.2 \cdot V_{No}(t)
\end{equation*}

where each control function is multiplied by the corresponding susceptible population.

We next set the mean number of people infected by an infectious person, namely the R0 number. Our model includes 4 age-specific R0 numbers associated with the $\beta$ parameters in our model, according to Equation (\ref{betar0}). In reality, the R0 number greatly fluctuates depending on individual behaviours and government restrictions. For that reason, we explore three different choices for the set of age-specific R0 numbers. These choices are listed in Table \ref{r0choices} along with references to the respective plots.

\begin{table}[ht]
	\centering
	\begin{tabular}{|c|c|c|c|c||c|c|} 
		\hline
		$R0_{OO}$&$R0_{YY}$&$R0_{OY}$&$R0_{YO}$&Average $R0$ &Plot before peak&Plot after peak\\
		\hline
		1&1&0.3&0.3&1.3&Figure \ref{1_and_0.3before}&Figure \ref{1_and_0.3after}\\
		1&1&0.9&0.9&1.9&Figure \ref{1_and_0.9before}&Figure \ref{1_and_0.9after}\\
		1.2&1.2&0.3&0.3&1.5&Figure \ref{1.2_and_0.3before}&Figure \ref{1.2_and_0.3after}\\
		\hline
	\end{tabular} 
	\caption{\label{r0choices}Chosen values for the $R0$ numbers (first four columns). The fifth column corresponds to the weighted average of the age-specific R0 numbers. The sixth and seventh column refer to the corresponding plots with the comparative simulation results, using Irish data before and after the peak in the number of infections respectively.}
\end{table}

\begin{figure}[!ht] 
	\begin{subfigure}{.5\textwidth}
		\includegraphics[width=\linewidth,height=4cm]{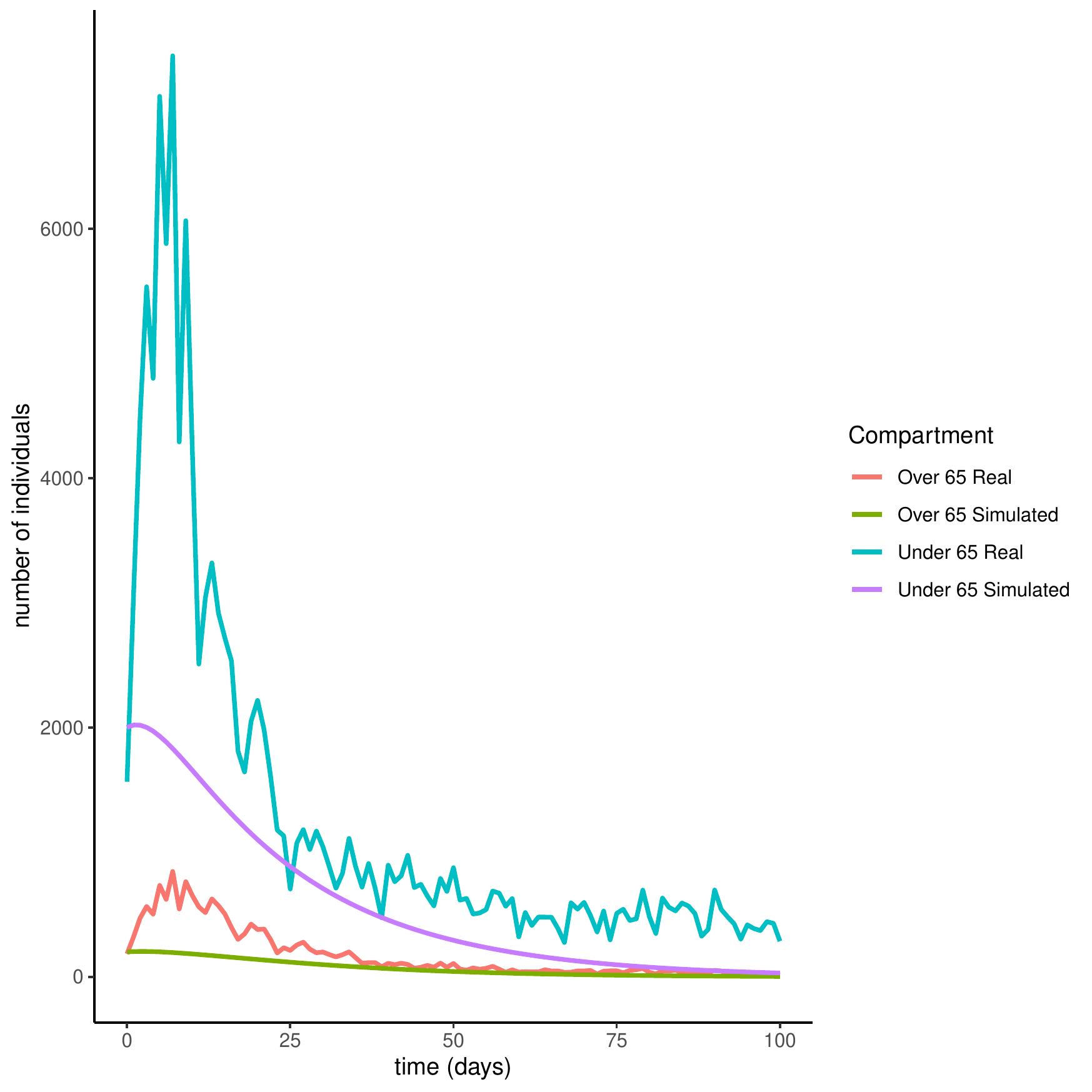}
		\caption{Number of infections}
	\end{subfigure}%
	\begin{subfigure}{.5\textwidth}
		\includegraphics[width=\linewidth,height=4cm]{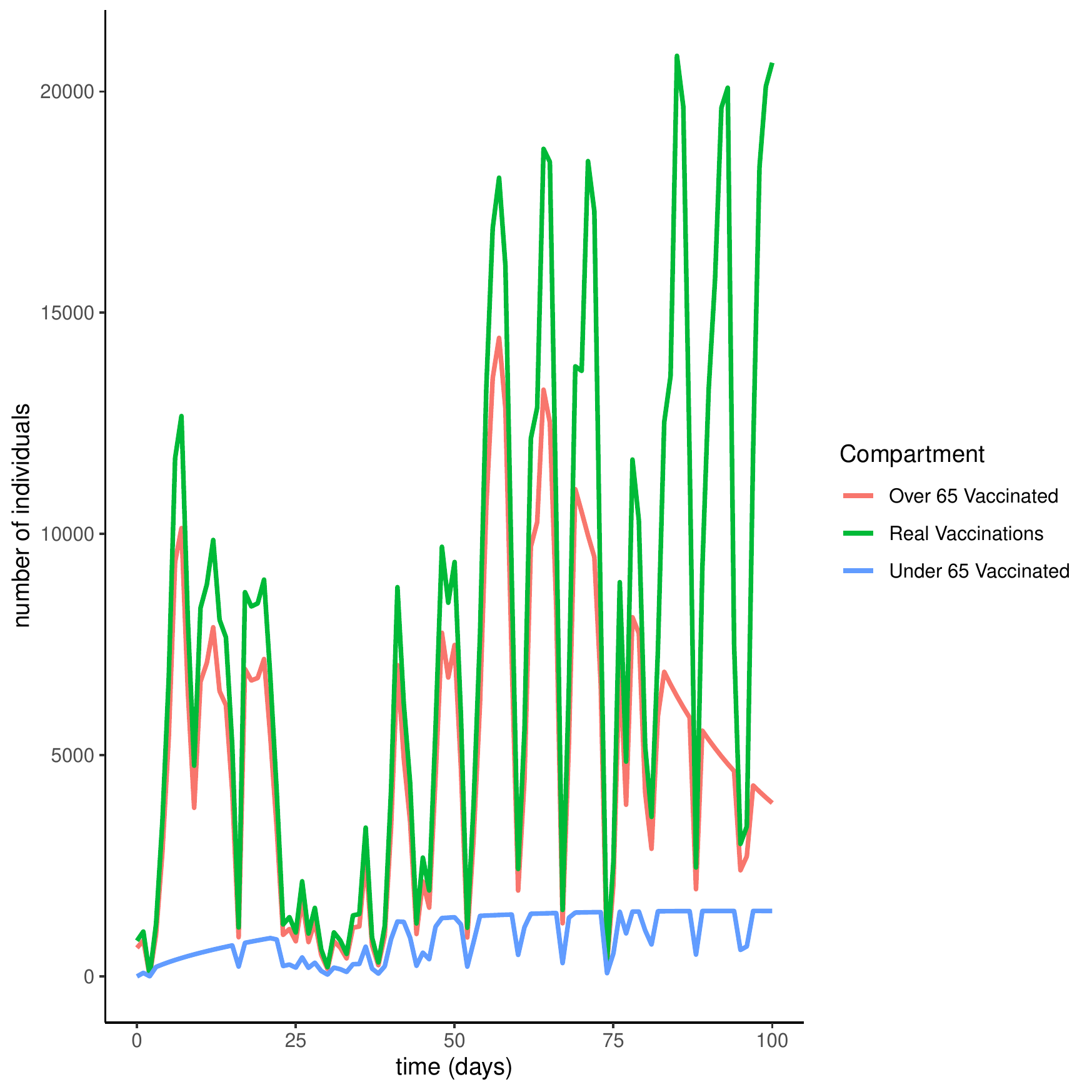}
		\caption{Number of vaccinations}
	\end{subfigure}
	\caption{\label{1_and_0.3before}Comparison of the infection and vaccination numbers for the Republic of Ireland starting from January 
		1\textsuperscript{st} 2021. The R0 numbers used are equal to 1 within each age group and 0.3 between age groups.}
\end{figure}

\begin{figure}[!ht] 
	\begin{subfigure}{.5\textwidth}
		\includegraphics[width=\linewidth,height=4cm]{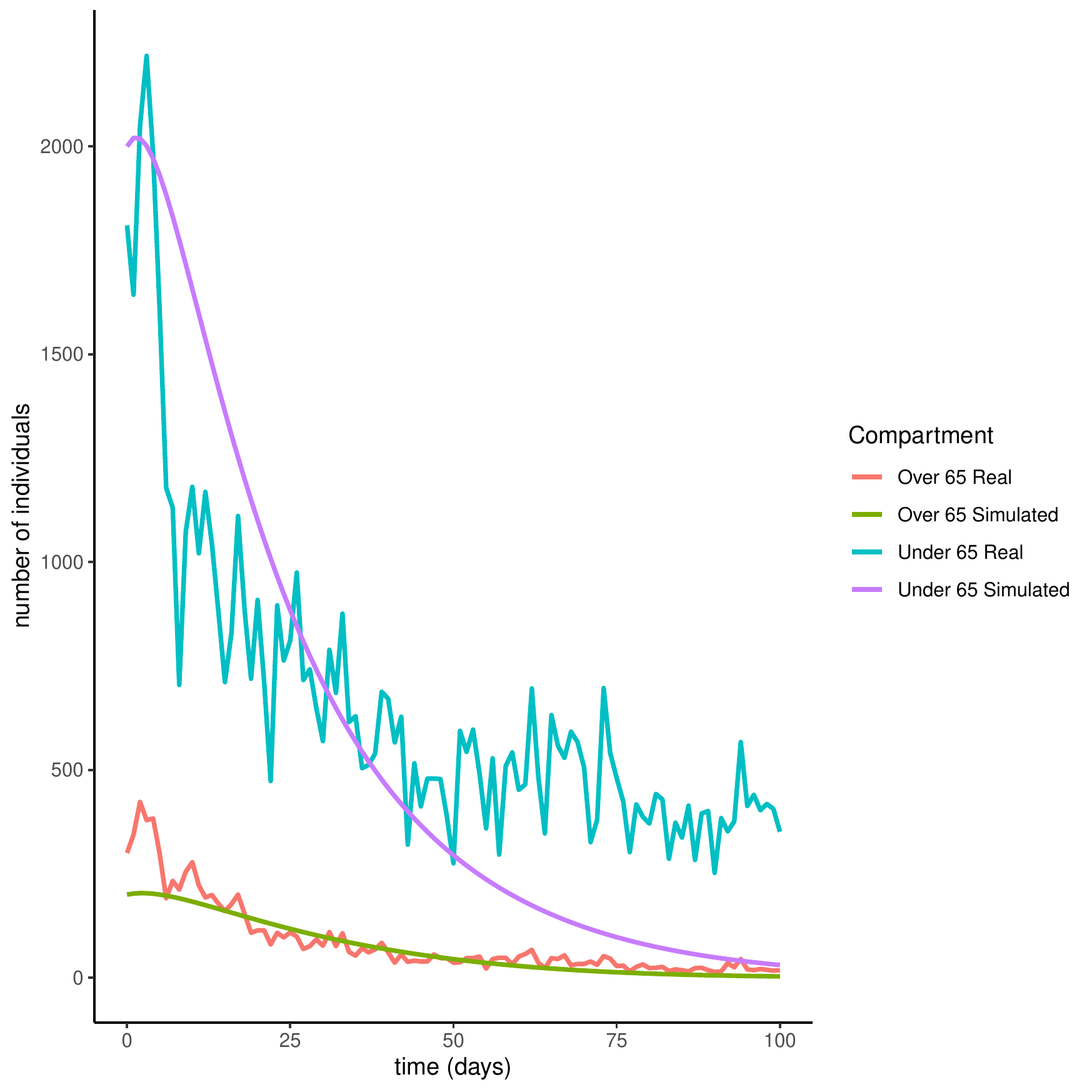}
		\caption{Infected individuals}
	\end{subfigure}%
	\begin{subfigure}{.5\textwidth}
		\includegraphics[width=\linewidth,height=4cm]{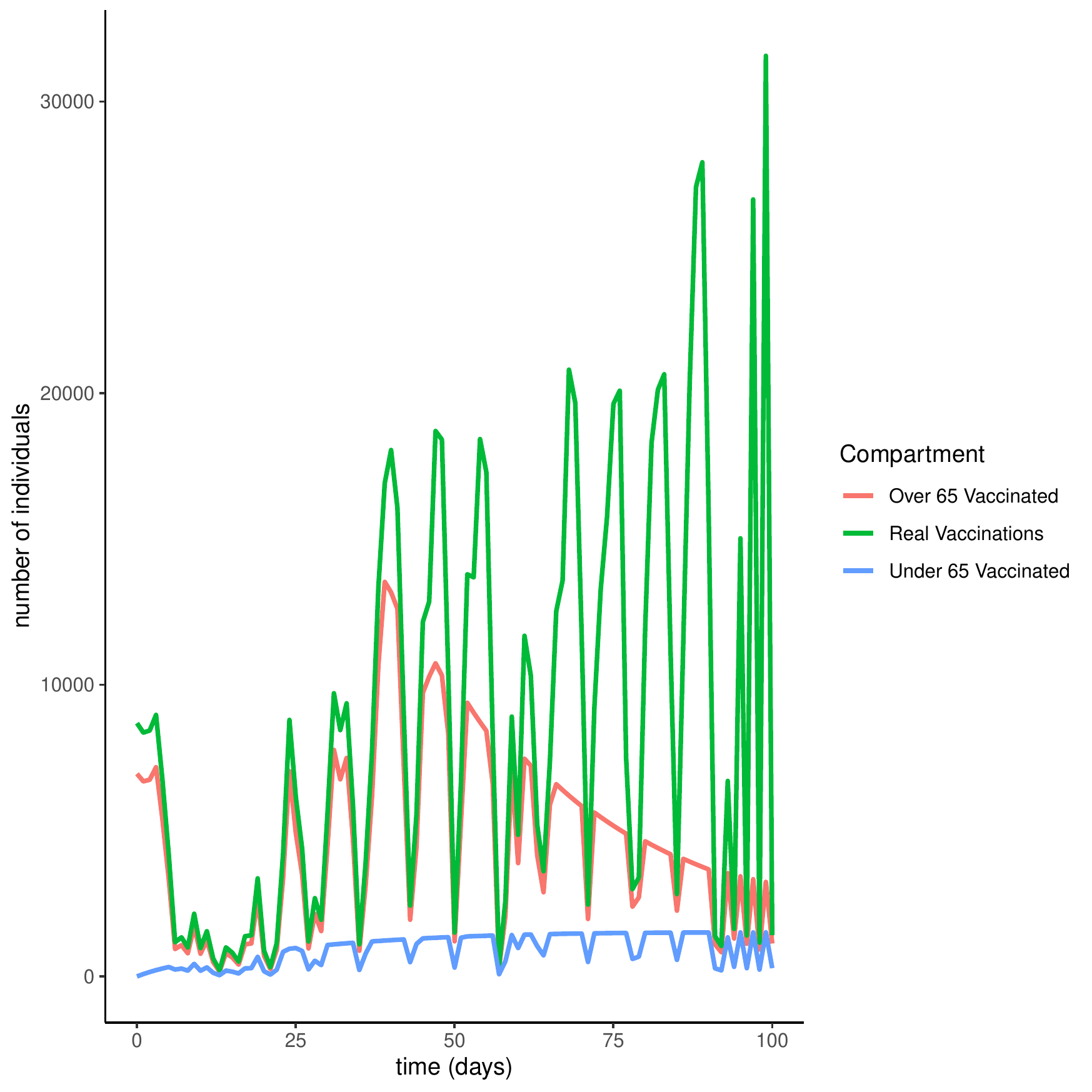}
		\caption{Vaccines administered}
	\end{subfigure}
	\caption{\label{1_and_0.3after}Comparison of the infection and vaccination numbers for the Republic of Ireland starting from January
		18\textsuperscript{th} 2021.The R0 numbers used are equal to 1 within each age group and 0.3 between age groups.}
\end{figure}

\begin{figure}[!ht] 
	\begin{subfigure}{.5\textwidth}
		\includegraphics[width=\linewidth,height=4cm]{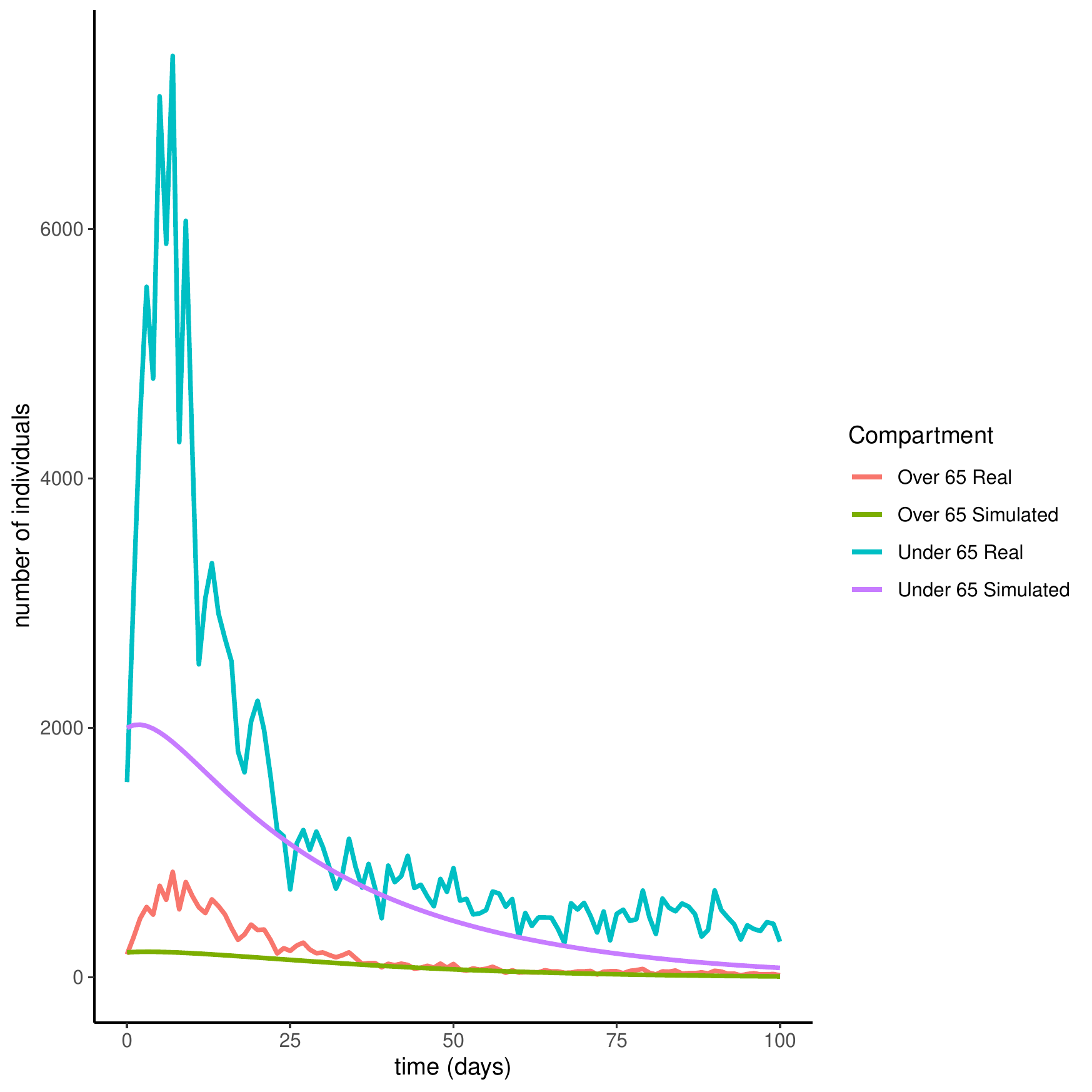}
		\caption{Number of infections}
	\end{subfigure}%
	\begin{subfigure}{.5\textwidth}
		\includegraphics[width=\linewidth,height=4cm]{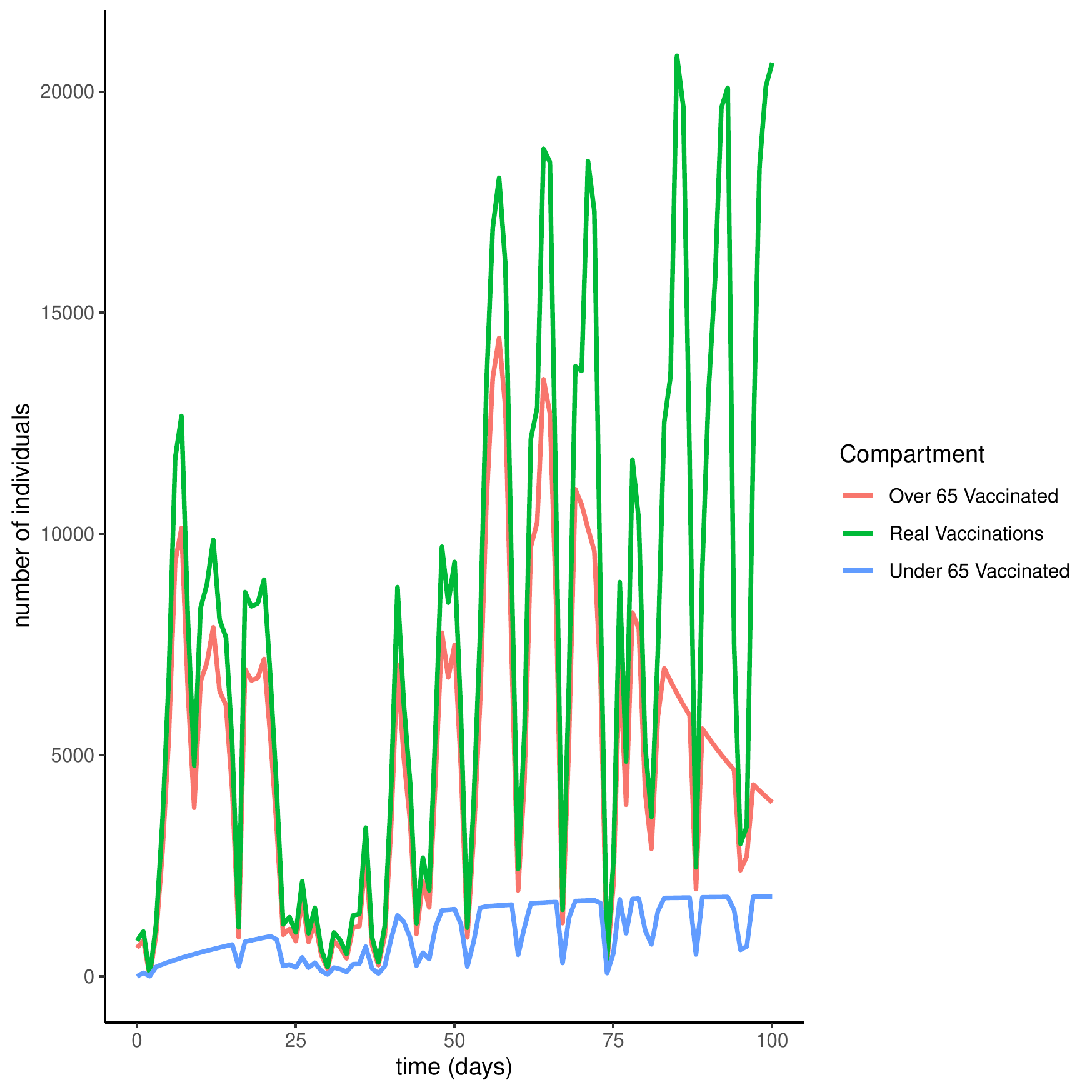}
		\caption{Number of vaccinations}
	\end{subfigure}
	\caption{\label{1.2_and_0.3before}Comparison of the infection and vaccination numbers for the Republic of Ireland starting from January 
		1\textsuperscript{st} 2021. The R0 numbers used are equal to 1.2 within each age group and 0.3 between age groups.}
\end{figure}

\begin{figure}[!ht] 
	\begin{subfigure}{.5\textwidth}
		\includegraphics[width=\linewidth,height=4cm]{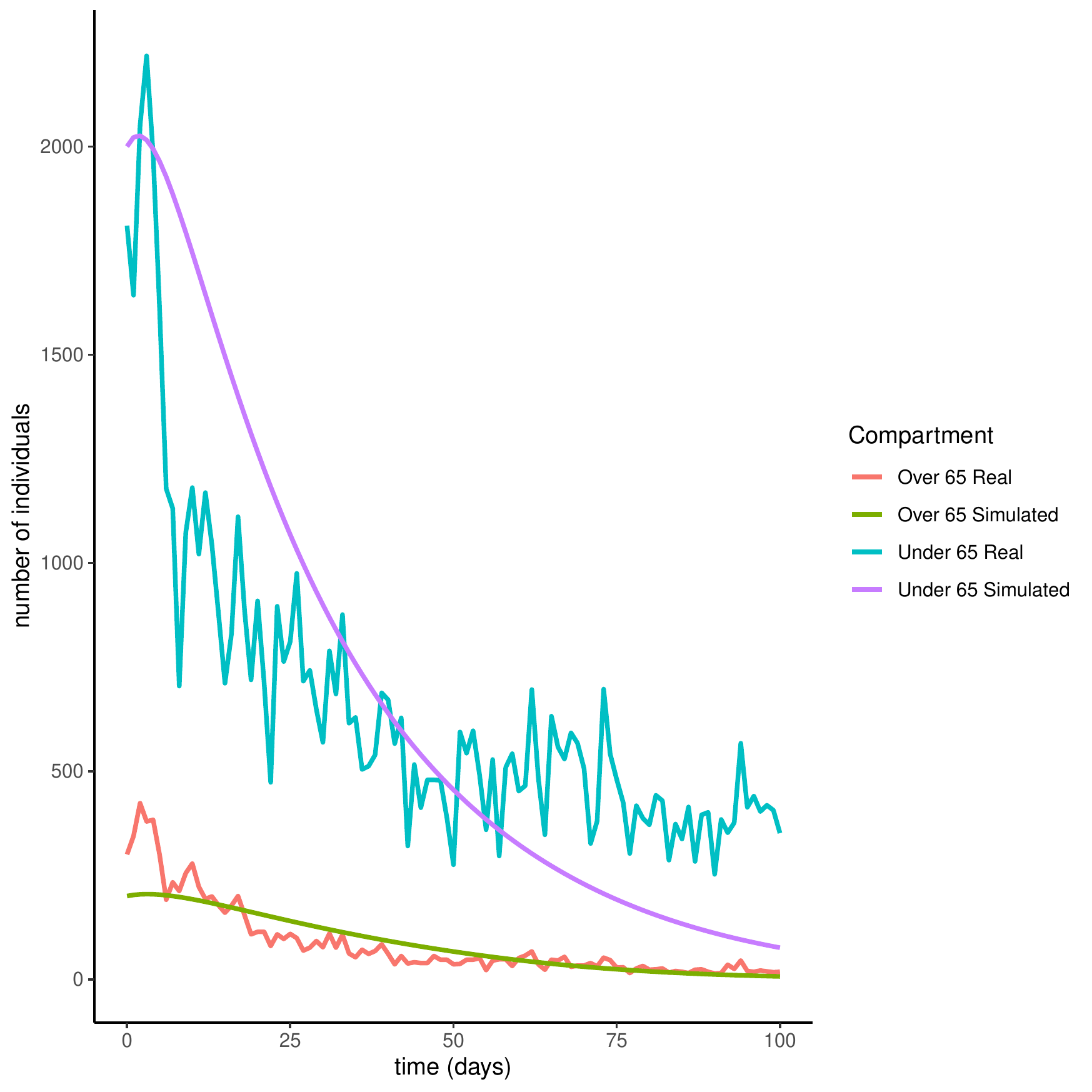}
		\caption{Infected individuals}
	\end{subfigure}%
	\begin{subfigure}{.5\textwidth}
		\includegraphics[width=\linewidth,height=4cm]{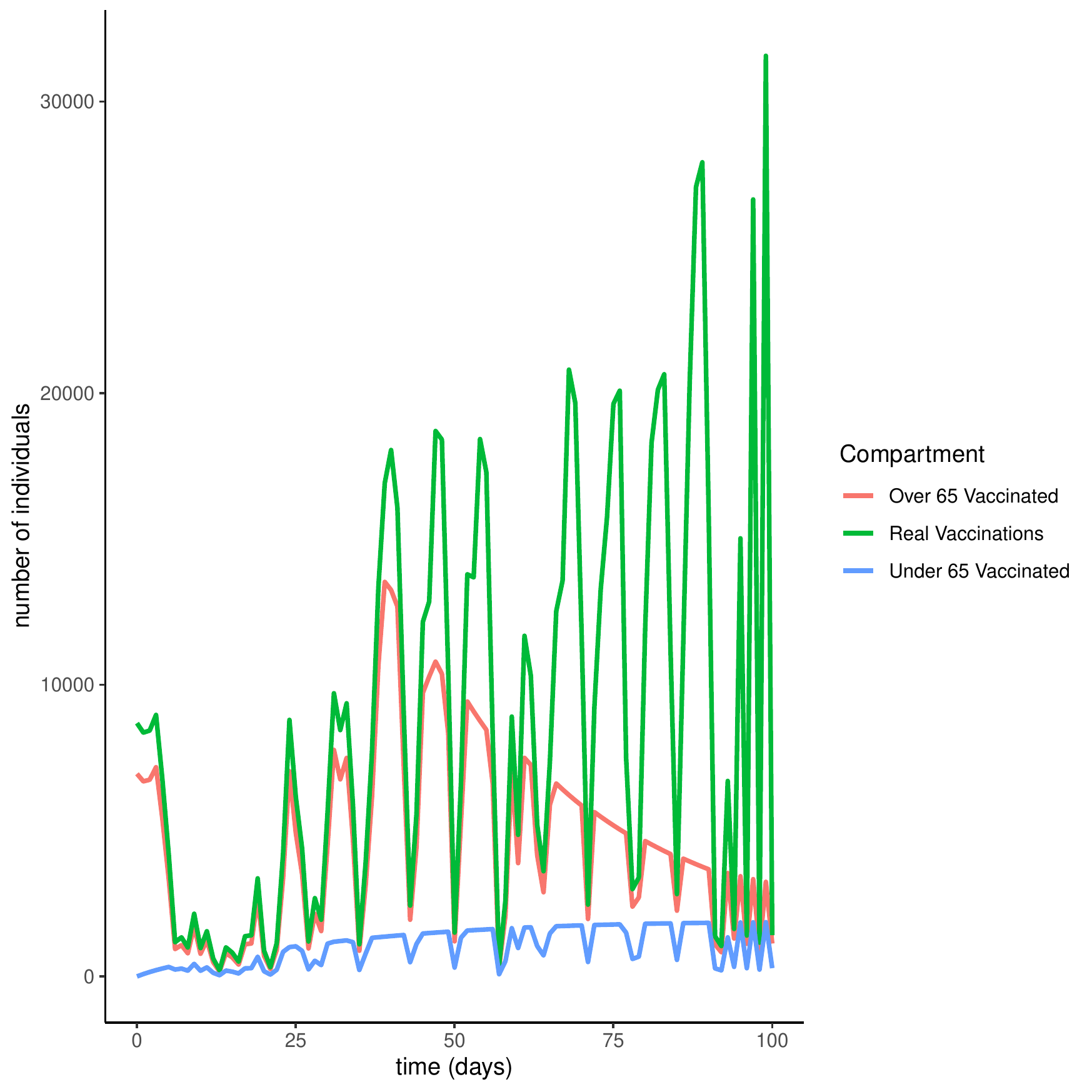}
		\caption{Vaccines administered}
	\end{subfigure}
	\caption{\label{1.2_and_0.3after}Comparison of the infection and vaccination numbers for the Republic of Ireland starting from January
		18\textsuperscript{th} 2021.The R0 numbers used are equal to 1.2 within each age group and 0.3 between age groups.}
\end{figure}

\begin{figure}[!ht] 
	\begin{subfigure}{.5\textwidth}
		\includegraphics[width=\linewidth,height=4cm]{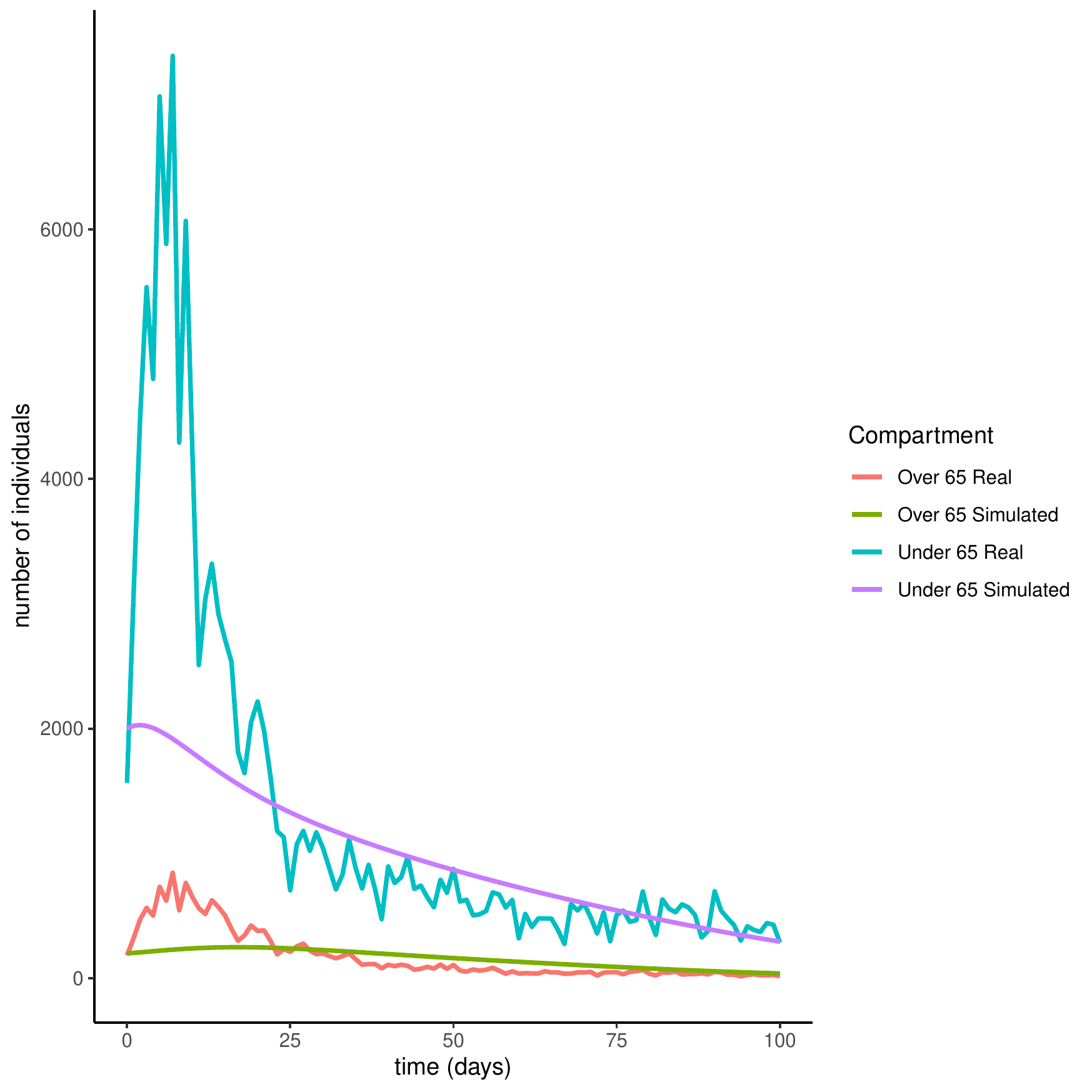}
		\caption{Number of infections}
	\end{subfigure}%
	\begin{subfigure}{.5\textwidth}
		\includegraphics[width=\linewidth,height=4cm]{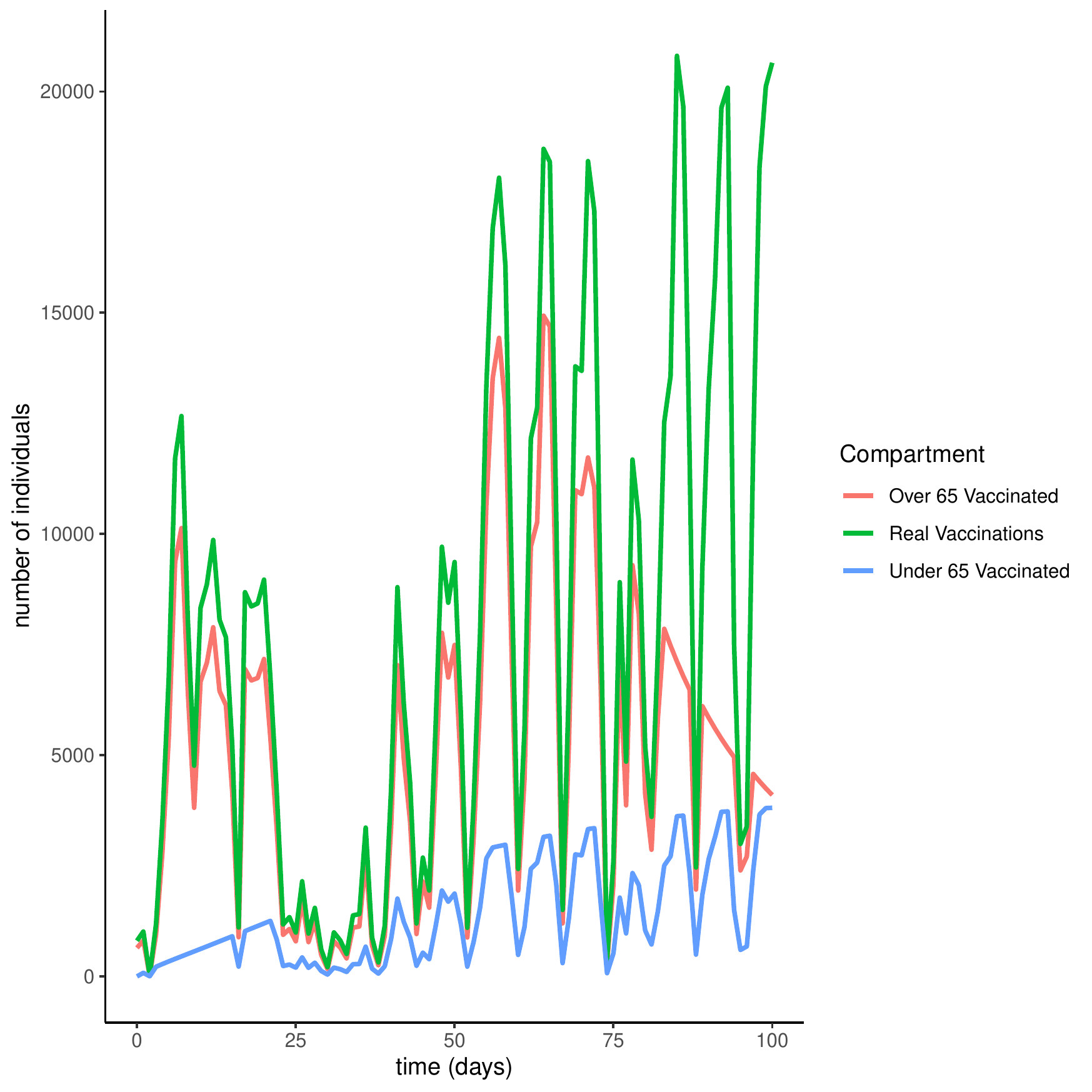}
		\caption{Number of vaccinations}
	\end{subfigure}
	\caption{\label{1_and_0.9before}Comparison of the infection and vaccination numbers for the Republic of Ireland starting from January 
		1\textsuperscript{st} 2021. The R0 numbers used are equal to 1 within each age group and 0.9 between age groups.}
\end{figure}

\begin{figure}[!ht] 
	\begin{subfigure}{.5\textwidth}
		\includegraphics[width=\linewidth,height=4cm]{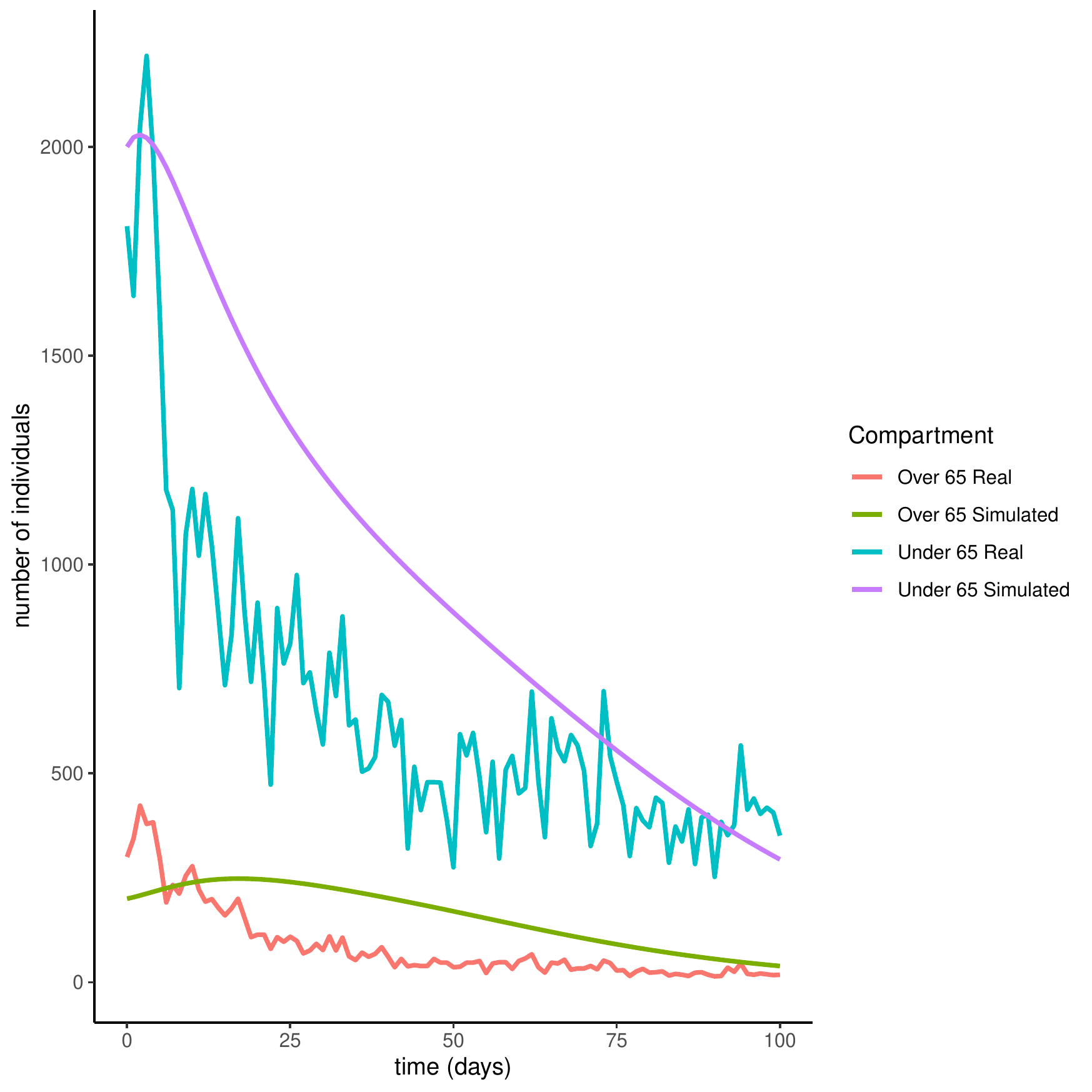}
		\caption{Infected individuals}
	\end{subfigure}%
	\begin{subfigure}{.5\textwidth}
		\includegraphics[width=\linewidth,height=4cm]{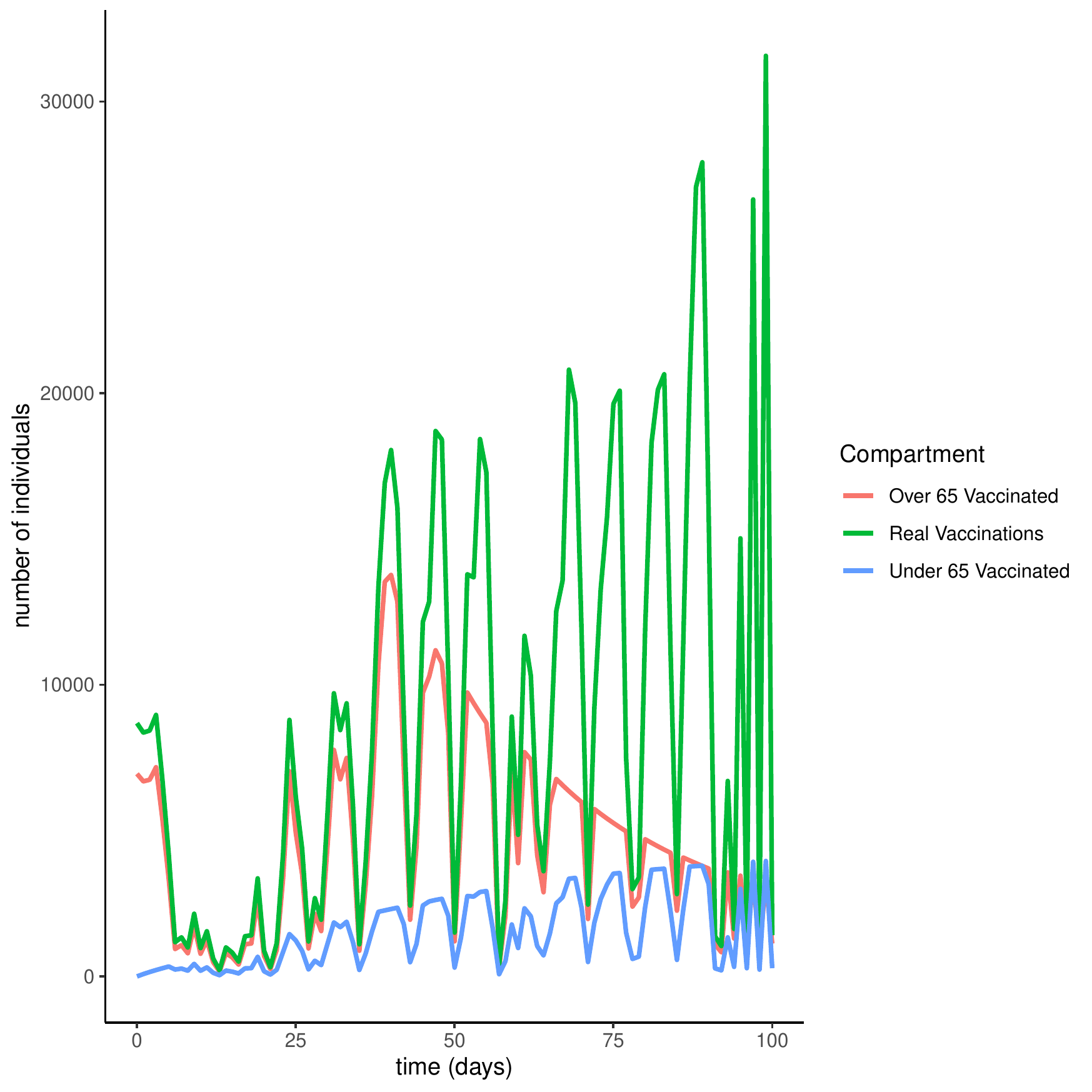}
		\caption{Vaccines administered}
	\end{subfigure}
	\caption{\label{1_and_0.9after}Comparison of the infection and vaccination numbers for the Republic of Ireland starting from January
		18\textsuperscript{th} 2021.The R0 numbers used are equal to 1 within each age group and 0.9 between age groups.}
\end{figure}

From Figures \ref{1_and_0.3before} to \ref{1.2_and_0.3after} it is evident that the optimal vaccination strategy obtained from our model is successful in quickly reducing the number of infections even in the case of an average R0 number of $1.9$. This value for R0 is unrealistic, especially for January 2021, when Ireland was under strict lockdown measures. However, in this extreme case, while the total number of infections is high compared to the real numbers, as expected, the decline of the infection curve is steep, as can be seen in Figure \ref{1_and_0.9after}. This implies that even if the country was not under total lockdown during these days, the optimal vaccination strategy would be able to bring the infections under control within a short amount of time. 

The successful control of infections for all choices of the $R0$ numbers can be attributed to the distribution of the vaccines to both age groups. The policy adopted by the Irish government mainly entailed administering vaccines to the over 65 population as are the higher risk group. However, it is evident from our simulation that a policy where a portion of the vaccines are additionally made available to the under 65 group can be more beneficial. 

It is interesting to note that in every simulation, the total number of daily vaccinations resulting from the optimisation are often less than the real vaccination numbers, i.e. the daily maximum number of available vaccines. This is a very encouraging observation, considering the implication that less vaccines can be used to achieve comparable or even better results. It is a particularly significant result if the optimisation is targeted towards countries that have gained limited access to vaccines. 

\section{Discussion} \label{section discussion}
Since the end of February 2020, the Republic of Ireland has been affected by the COVID-19 pandemic, like most countries in the world. For the greater part of the year no effective treatment or vaccine was present, leaving the government with one tool to attempt to flatten the curve of infection: the application of various levels of restrictions, depending on the occurrence of new daily cases. In light of the vaccines that became available near the end of 2020, a need arose for efficient vaccination plans based on the availability of the vaccines and the risk levels of different groups of the population.

In this work, we introduced a novel compartmental model which more adequately describes the dynamics of the virus compared to the approaches taken thus far in the literature. We added compartments representing the stages of infection and vaccination, and assumed two different age groups sharing the same model structure. Using this model as a baseline to describe the evolution of the pandemic in the Republic of Ireland, we applied Optimal Control methods to obtain a suggested vaccination strategy that minimises the number of infections under certain restrictions such as the available vaccines. We simulated the evolution of the pandemic with and without the vaccination strategy in place for two different scenarios, indicative of strict and loose restrictions respectively. The main conclusions we drew were:
\begin{enumerate}
	\item The optimal control strategy in both scenarios successfully flattens the curve of infectious individuals, ensuring that less people will get infected by the virus due to protection from the vaccine and, most importantly, that less people will be simultaneously infected thus avoiding exhausting the national healthcare system. 
	\item In the case of low transmission rates, the vaccination not only flattens the curve but also significantly reduces the total number of infections, specifically from 80\% to 11.3\% for the over 65 group and from 79\% to 5.38\% for the under 65 group. However, the pandemic lasts a long time, meaning that a big part of the population remains in the susceptible compartment, not yet vaccinated but also not infected due to the low transmission rates, as opposed to being moved to the protected or recovered compartments. 
	\item When the transmission rates are high however, the vaccination reduces the number of infections by a much smaller percentage, namely from 99.995\% to 85.29\% for the over 65 group and from 99.993\% to 99.1\% for the under 65 group. However it still succeeds in flattening the curve and ensuring less simultaneous infections, hence preventing a healthcare crisis. The pandemic in this case ends in a very short period of time with the biggest part of the population having been infected and eventually recovered.
	\item For both scenarios, the optimal strategy suggests focusing on vaccinating the older population in higher percentages first, while simultaneously vaccinating part of the younger population. This is a different approach to the one taken by many states, including Ireland, where the older population gets exclusively vaccinated first. 
\end{enumerate}
All of the above lead us to conclude that an approach involving a combination of strict and loose measures would be ideal, while the vaccination programme is taking place. That would ensure that the infection doesn't spread to such a high percentage of the population while the restrictions are not as severe so as to make the situation unbearable for a certain length of time.

To verify our findings, we compared the results of our simulations to the real infection numbers that occurred during the beginning of 2021 in Ireland, when the vaccination process was starting across the country. We used three different sets of transmission numbers, each of them representative of different levels of restrictions. Our method was successful in bringing the infection numbers under control within a short amount of time and with a need for less resources than what was used. This result is particularly important as it stresses the importance of correct allocation of the available resources over the quantity of those resources. This means that even countries with limited access to vaccines are able to effectively reduce the number of infections.

While our approach led us to a lot of interesting conclusions, there are certain drawbacks. Our model assumes that all of the parameters remain constant for the duration of the simulations which is not a realistic assumption. For example the transmission rates, which were the source of the vastly different results in our two simulations, do not remain constant for long periods of time. Furthermore, our model is a deterministic model which means that there is no accounting for any uncertainty. However it would be relatively simple to propagate uncertainty through the model by repeatedly running it with parameter values sourced from an appropriate probability distribution.

In order to avoid the drawbacks described above, a stochastic model would have to be introduced or a model with variable parameters. A study on a time-varying stochastic SEIR model for the control of the Ebola virus can be found in \cite{lekone2006statistical}. In addition, as we mentioned earlier, it is interesting to explore alternative expressions for the objective functional to be minimised. For example, different weights to the control functions for each age group can be introduced, modelling the cost of bringing the vaccine to different populations, or weight factors relating to the mortality rates of each group, to the numbers of infectious individuals. Moreover, the introduction of more compartments could enrich the model further, for instance a compartment for the hospitalised individuals or one for the deceased. Finally, notable extensions would be the introduction of more age groups and high risk groups, as well as multiple vaccines and their effectiveness. We leave these topics as future research topics which we pursue elsewhere.

\section*{Acknowledgments}

This publication has emanated from research supported in part by a research grant from Science Foundation Ireland (SFI) under Grant Number 16/RC/3872 and is co-funded under the European Regional Development Fund. This work was also supported by a Science Foundation Ireland COVID-19 Rapid Response Grant number 20/COV/0081. Andrew Parnell’s work was additionally supported by: a Science Foundation Ireland Career Development Award (17/CDA/4695); an investigator award (16/IA/4520); a Marine Research Programme funded by the Irish Government, co-financed by the European Regional Development Fund (Grant-Aid Agreement No. PBA/CC/18/01); European Union’s Horizon 2020 research and innovation programme InnoVar under grant agreement No 818144; SFI Centre for Research Training in Foundations of Data Science 18CRT/6049, and SFI Research Centre award 12/RC/2289\_P2.

\newpage
\appendix 
\section{Table of notations}\label{appendix_notations}
\begin{center}
	\begin{tabular}{|c|c|}
		\hline
		\textbf{Model states} & \textbf{Symbol} \\
		\hline
		Susceptible not yet vaccinated (o65 and y65) & $S_O, S_Y$ \\
		Received vaccine, waiting for it to take effect (o65 and y65)& $V_O, V_Y$\\
		Received vaccine but was not effective (o65 and y65)& $N_O,N_Y$ \\
		Susceptible, refusing or unable to receive vaccine (o65 and y65)& $U_O,U_Y$ \\
		Exposed (o65 and y65)& $E_O,E_Y$\\
		Infectious (o65 and y65)& $I_O,I_Y$ \\
		Recovered (o65 and y65)& $R_O,R_Y$ \\
		Protected from vaccine (o65 and y65)& $P_O,P_Y$\\
		Total number of people in age group (o65 and y65) & $T_O,T_Y$\\
		\hline
		\textbf{Model Parameters} & \textbf{Symbol} \\
		\hline
		Rate at which an o65 person infects an o65 person & $\beta_{OO}$ \\
		Average number of o65 people infected by an o65 person & $R0_{OO}$\\
		Rate at which a y65 person infects an o65 person & $\beta_{YO}$\\
		Average number of o65 people infected by a y65 person & $R0_{YO}$\\
		Rate at which an o65 person infects a y65 person & $\beta_{OY}$\\
		Average number of y65 people infected by an o65 person & $R0_{OY}$\\
		Rate at which a y65 person infects a y65 person & $\beta_{YY}$\\
		Average number of y65 people infected by a y65 person & $R0_{YY}$\\
		Rate at which exposed becomes infected &$\gamma_E$ \\
		Rate at which infected becomes recovered &$\gamma_I$ \\
		Rate at which vaccinated becomes protected &$\gamma_V$ \\
		Vaccine effectiveness & $\alpha_V$\\
		Percentage of over 65s refusing the vaccine & $r_O$\\
		Percentage of under 65s refusing the vaccine & $r_Y$\\
		\hline 
		\textbf{Control functions}& \textbf{Symbol}\\
		\hline
		Percentage of over 65s to get vaccinated at time $t$ & $u_O(t)$\\ 
		Percentage of under 65s to get vaccinated at time $t$ & $u_Y(t)$\\
		\hline
		\textbf{Optimisation Functions and Parameters}& \textbf{Symbol}\\
		\hline
		Hamiltonian function &H \\
		Cost function & $\mathcal{F}$\\
		Age specific weight constants & $W_O,W_Y$\\
		\hline
\end{tabular}\end{center}
\label{table of notations}

\newpage
\section{Existence of an optimal solution}\label{appendix_conditions}

The state equations \eqref{stateequations} can be written in compact form as
\begin{equation}
g(t,x,u):=\frac{dx(t)}{dt}
\end{equation}
where $x$ denotes the variable vector $x= \left(\right.S_O,V_O,N_O  ,U_O,E_O,I_O,R_O,P_O,S_Y,V_Y,N_Y,U_Y,$ $E_Y,I_Y,R_Y ,P_Y\left.\right) $. Firstly, we assume $g$ satisfies the following:

\begin{enumerate}[(a)]
	\item \label{assum_1}$\vert g(t,x,u) \vert \leq C_1 (1+\vert x\vert + \vert u \vert )$
	\item \label{assum_2}$\vert g(t,x,u)- g(t,x',u) \vert \leq C_2 \vert x-x'\vert (1+\vert u \vert )$
\end{enumerate}
for some positive constants $C_1, C_2$ and for all $t,x,x'$ and $u\in \Omega$. As stated in \cite{Fleming1975}, when $g$ is of the class $C^1$, meaning that it is differentiable with continuous first order partial derivatives, then the above conditions are satisfied by applying suitable bounds on the partial derivatives of $g$. It is possible to show that such bounds can be applied in the case of our epidemic model, from the boundedness of the control and state variables ensuring that 
Assumptions (\ref{assum_1})-(\ref{assum_2}) are true.

Consider a performance index of the form: 
\begin{equation}
\mathcal{F}(u(t)) = \int_0^T f(t,x(t),u(t))dt
\end{equation}
The following theorem, provides the sufficient conditions for the existence of a solution to the optimal control problem.
\begin{theorem} \cite{Fleming1975}
	Suppose that assumptions \ref{assum_1}-\ref{assum_2} hold, that $f$ is continuous and moreover that:
	\begin{enumerate}
		\item \label{cond_1} The set of solutions to the system equations (\ref{stateequations}) with corresponding control functions in $\Omega$ is nonempty.
		\item \label{cond_2}  $\Omega$ is convex.
		\item  \label{cond_3} $f$ is convex on $\Omega$ and $g$ can be written as $g(t,x,u)=a(t,x)+b(t,x)u$.
	\end{enumerate}
	Then there exists an optimal control $u^*(t)$ with corresponding optimal trajectory $x^*(t)$ minimizing $\mathcal{F}(u(t))$.
\end{theorem}
We now prove the conditions of the theorem are satisfied for our model. As noted earlier, conditions (\ref{assum_1})-(\ref{assum_2}) are satisfied. The performance index integrand is 
\begin{equation}
f(t,x(t),u(t))= I_O(t)+ I_Y(t) + \frac{W_O}{2}u_O^2(t) + \frac{W_Y}{2}u_Y^2(t) 
\end{equation}
which is continuous, and $\Omega = \left\lbrace u(t) \in L^2(O,T)^2 \Vert a \leq u_O(t),u_Y(t)\leq b, t \in [0,T]\right\rbrace$. We now examine conditions \ref{cond_1}-\ref{cond_3}.
\begin{enumerate}
	\item From the \textit{Picard-Lindelöf} theorem \cite{coddington1955theory}, if the system equations are continuous and Lipschitz and the solutions to the system equations are bounded, then there is a unique solution corresponding to every admissible control in $\Omega$. 
	It is trivial to verify the continuity of the state equations \eqref{stateequations}. Additionally, every element of the variable vector $x(t)$ is bounded by the corresponding age-group totals $T_O, T_Y$. In other words, $x(t)\in [O,\max\{T_O,T_Y\}]^{16}$ which is a compact and convex set. Finally, the Lipschitz property can be verified with the use of the following lemma from \cite{birkhoff1978ordinary}.
	\begin{lemma}
		If the vector function $\mathbf{X}(\mathbf{x},t)$ is of class $C^1$ in a compact convex domain D, then it satisfies the Lipschitz condition there.
	\end{lemma}
	
	\item From the \textit{Heine-Borel} theorem, it follows that $\Omega$ is convex because as a closed and bounded set in $\mathbb{R}^2$ \cite{rudin1976principles}.
	\item The state equations are linearly dependent on the control functions $u_O(t)$ and $u_Y(t)$ so they can be written in the form $g(t,x,u)=a(t,x)+b(t,x)u$. The convexity of the integrand $f$ in the objective functional, follows from the fact that it is quadratic in the control functions.
\end{enumerate}

\bibliography{mybibfile}
\bibliographystyle{plain}

\end{document}